\documentclass[smallextended]{svjour3}
\makeatletter 
\let\cl@chapter\undefined
\makeatletter
\usepackage[english]{babel}

\usepackage[letterpaper,top=2cm,bottom=2cm,left=3cm,right=3cm,marginparwidth=1.75cm]{geometry}

\usepackage{natbib}
\usepackage{xparse}
\usepackage{amsmath}
\usepackage{amssymb}
\usepackage{mathtools}
\usepackage{graphicx}
\usepackage{subcaption}
\usepackage[colorlinks=true, allcolors=blue]{hyperref}
\usepackage{pgfplots}
\usepackage{import}
\usepackage{color}
\usepackage{xcolor}
\usepackage{csquotes}
\usepackage{multirow}
\usepackage{todonotes}
\usepackage{varioref}
\usepackage[capitalise,noabbrev]{cleveref}
\usepackage[title]{appendix}
\usepackage{float}
\usepackage{booktabs}

\pgfplotsset{compat=1.16}

\title{An analysis of optimization problems involving ReLU neural networks}
\author{Christoph Plate \and 
            Mirko Hahn  \and
            Alexander Klimek \and
            Caroline Ganzer \and
            Kai Sundmacher \and
            Sebastian Sager}
\institute{Christoph Plate, Mirko Hahn, Kai Sundmacher, Sebastian Sager \at
         Otto von Guericke University Magdeburg, Magdeburg, Germany \\
         \email{\{christoph.plate, mirhahn, kai.sundmacher, sager\}@ovgu.de}
         \and
         Christoph Plate, Alexander Klimek, Caroline Ganzer, Kai Sundmacher, Sebastian Sager \at
         Max Planck Institute for Dynamics of Complex Technical Systems, Magdeburg, Germany \\
         \email{\{plate, klimek, cganzer, sundmacher, sager\}@mpi-magdeburg.mpg.de}
         \and
         Corresponding author: Christoph Plate
}

\date{Received: date / Accepted: date}
\newcommand{\myReLU}[1]{\textrm{ReLU}\left(#1\right)}
\newcommand{\myClippedReLU}[2][M]{\textrm{ReLU}_{#1}\left(#2\right)}

\begin{document}
\maketitle

\begin{abstract}
Solving mixed-integer optimization problems with embedded neural networks with ReLU activation functions is challenging.
Big-M coefficients that arise in relaxing binary decisions related to these functions grow exponentially with the number of layers. 
We survey and propose different approaches to analyze and improve the run time behavior of mixed-integer programming solvers in this context.
Among them are clipped variants and regularization techniques applied during training as well as optimization-based bound tightening and a novel scaling for given ReLU networks. 
We numerically compare these approaches for three benchmark problems from the literature.
We use the number of linear regions, the percentage of stable neurons, and overall computational effort as indicators. 
As a major takeaway we observe and quantify a trade-off between the often desired redundancy of neural network models versus the computational costs for solving related optimization problems.
\end{abstract}

\section{Introduction}\label{sec:introduction}

Artificial neural networks (ANNs) are a popular tool for approximating functions from data and have been used in various applications, ranging from modeling and control of batch reactors \citep{Mujtaba2006}, the optimization of cancer treatments \citep{Bertsimas2016} and chemical reactions \citep{Fernandes2006} to the approximation of solutions to complex optimization problems \citep{Bertsimas2022}. 
Formally, a feed-forward ANN consists of $J \in \mathbb N$ consecutive layers. Each layer uses an affine-linear transformation of the input with a subsequent, element-wise application of a nonlinear activation function $s: \mathbb R \mapsto \mathbb R$, i.e.,
\begin{align}\label{eq:ANN}
	x^{(j)} = s^{(j)}\left( W^{(j)} x^{(j-1)} + b^{(j)}  \right), \quad j \in [J],
\end{align}
with $x^{(0)} = x \in \mathbb{R}^{n_x}$ as the input to the neural network and $W^{(j)} \in \mathbb{R}^{n_j\times n_{j-1}}, b^{(j)} \in \mathbb{R}^{n_j}$ denoting the weights and biases of layer $j$, respectively.
In this paper, we assume $s$ to be the ReLU activation
\begin{align}\label{eq:relu}
    \myReLU{x} = \max\{0,x\},
\end{align}
on all hidden layers and the identity on the last layer. 
A recent survey shows that more than 400 different activation functions have been suggested in the literature \citep{kunc2024three}. While in some contexts the usage of other, possibly continuously differentiable, activation functions may be recommendable, ReLU activation variants continue to play a major role. Especially in the context of mixed-integer optimization problems where some of the variables require a non-smooth, non-differentiable treatment, a study of ReLU activation is thus of major interest.
Equipped with the ReLU activation function, an ANN describes a piecewise affine-linear function $h \colon \mathbb{R}^{n_x} \mapsto \mathbb{R}^{n_h}$ with $h(x)=x^{(J)}$~\citep{Grigsby2022}. Under mild assumptions on their size and the chosen activation function, neural networks are universal approximators, i.e., they can approximate continuous functions on compact sets to arbitrary precision \citep{Cybenko1989,Hornik1991}. Therefore, they are often used for approximating functional relationships in cases where the underlying function is unknown, hard to model otherwise, or in general expensive to simulate, but where data is available to train the network on. See, e.g., \citet{Misener2023} for a survey on surrogate modeling in process applications. Further, for classification tasks, such as image classification, ReLU ANNs and convolutional neural networks in particular play an important role \citep{Krizhevsky2012}.
Due to their widespread use, embedding trained universal approximators in optimization problems and the efficient solution of these problems is of special interest and currently an active field of research \citep{Schweidtmann2019a,Tong2024}. The applications in which these optimization problems appear are manifold. 
In general, the combination of first-principle models with data-driven surrogates has many advantages \citep{camps2023discovering}. We thus expect the optimization of mathematical models involving neural networks to play an important role in future engineering research, compare also \citet{Schweidtmann2021} for the case of chemical engineering.
We are interested in the setting where an ANN is embedded in a mixed-integer nonlinear optimization problem (MINLP), i.e.,
\begin{equation}\label{prob:embedded}
    \begin{alignedat}{3}
        \underset{x \in \mathbb{R}^{n_x}, w \in \mathcal{W}}\min&\ && f(y,w) \\
        \underset{\phantom{x \in \mathbb{R}^{n_x}, w \in \mathcal{W}}}{\textrm{s.t.}}&\ && \begin{aligned}[t]
            y & = h(x),\\
            0 & \leq g(y,w),
        \end{aligned}
    \end{alignedat}
\end{equation}
where $x \in \mathbb{R}^{n_x}$ and $y \in \mathbb{R}^{n_h}$ are input and output variables to the neural network, respectively, and $\mathcal{W}$ is a feasible set of $n_w$ additional mixed-integer variables $w$. Note that $w$ does not influence the neural network, but the optimization problem. The objective function $f \colon \mathbb{R}^{n_h} \times \mathcal{W} \mapsto \mathbb{R}$ and the constraint function $g \colon \mathbb{R}^{n_h} \times \mathcal{W} \mapsto \mathbb{R}^{n_c}$ typically depend on the output of the neural network and on the variables $w$. This general formulation is easily extendable towards the embedding of multiple neural networks, but also includes the simple case of minimizing the output of a single ANN in the absence of variables $w$ or additional constraints $g(\cdot)$. 
The task of verifying the reliability of neural networks can also be reduced to solving optimization problems of this type. Such verification problems arise because ANNs can be prone to adversarial attacks, i.e., situations in which minor perturbations in the input cause the network to produce incorrect outputs. Examples of such behavior can be found by solving the optimization problem 
\begin{equation}\label{prob:adversarial}
    \begin{aligned}
        \max_{\varepsilon \in \mathbb{R}^{n_x}}\ & h(x+\varepsilon)_k - h(x+\varepsilon)_i \\
        \textrm{s.t.}\ & \lvert \varepsilon \rvert \leq \delta,
    \end{aligned}
\end{equation}
where $\varepsilon$ is the perturbation which, when added to a current input $x$, changes the prediction from the correct label $i$ to the incorrect label $k$. The bound on the perturbation $\delta \in \mathbb{R}$ and its norm are hyperparameters of this problem and must be chosen on a problem-by-problem basis. See, e.g., \citet{Hein2017} for examples using the $\ell_2$ norm, or \citet{Tjeng2019} using the  $\ell_\infty$ norm. When ANNs are used in safety-critical applications, solving \eqref{prob:adversarial} is important to certify the robustness of the ANN. For instance, it can be used to verify that no adversarial example exists around specific inputs, e.g., those in the training set. Robustness against adversarial attacks can be proven if the optimal objective value of \eqref{prob:adversarial} is negative. See also \citet{Rossig2021} and the references therein for more information on verification problems. 

Embedding the ANNs into optimization problems entails introducing the necessary variables and modeling the equations \eqref{eq:ANN} for each neuron in the neural network. 
The first detailed study of the impact on the optimization problem and a comparison of different modeling approaches and algorithms was given in \citet{Joseph-Duran2014}. The application setting was different, because the functions $\max(0, x)$ modeled the overflow of sewage water. Yet, as an identical function to \eqref{eq:relu} was used, the mathematical formulation is identical. A major insight was that tailored constraint branching algorithms outperform standard mixed-integer modeling and continuously differentiable reformulations of \eqref{eq:relu}.
Nevertheless, for ReLU activations a big-M formulation is most widely used in the literature \citep{Fischetti2018,Xiao2019,Tjeng2019}. Big-M formulations are a standard technique to model constraints that can be activated or deactivated in mixed-integer linear programming (MILP). For a single neuron $i$ in layer $j$ $x^{(j)}_i = \text{ReLU}(W^{(j)}_i x^{j-1} + b_i)$, and assuming a bounded input $L^{(j)}_i \leq W^{(j)}_i x^{(j-1)} + b_i\leq U^{(j)}_i$, the big-M formulation reads
\begin{align}
    \begin{split}
        x^{(j)}_i & \geq  0, \\
        x^{(j)}_i & \geq W^{(j)}_i x^{(j-1)} + b_i,\\
        x^{(j)}_i & \leq  W^{(j)}_i x^{(j-1)} + b_i - L^{(j)}_i (1-z_i^{(j)}), \\
        x^{(j)}_i & \leq U^{(j)}_i z_i^{(j)}, \\
        z^{(j)}_i & \in \{0,1\},
    \end{split}
    \label{eq:bigM}
\end{align}
where $W^{(j)}_i$ denotes the $i$-th row of the weight matrix in layer $j$. Although this formulation is easily derived and implemented, the choice of the big-M coefficients $L^{(j)}_i$ and $U^{(j)}_i$ is crucial for practical performance. Larger values lead to weaker relaxations and can slow the convergence of MINLP solvers.
There is ongoing research to derive formulations with tighter relaxations or problem-specific cuts. An extended formulation proposed in \citet{Anderson2020}, which can be proven to yield the tightest possible relaxation for each neuron. This comes at the price of introducing additional continuous variables. However, the authors' own numerical studies find that the extended formulation does not offer significant performance improvements in optimization despite its theoretical advantages over the big-M formulation. A class of intermediate formulations between the big-M formulation and the extended formulation were proposed in \citet{Tsay2021,Kronqvist2024}. These formulations allow for a trade-off between dimension and relaxation tightness. In the extreme cases, they correspond exactly to the two formulations. The authors demonstrate with numerical experiments that their proposed partition-based formulation performs better in some application settings. For more information on ReLU ANNs and their MILP encodings we refer to the extensive survey \citet{Huchette2023a} and the references therein.

Several methods have been proposed to reduce the computational burden of solving optimization problems with embedded neural networks. These fall broadly into two categories. 

First, the ANN training can be adapted to yield ANNs with properties that ease the subsequent optimization. In \citet{Xiao2019}, the authors discuss regularization methods that can be applied during the training of the neural network which significantly speed up the solution of subsequent verification problems. Besides standard $\ell^1$ regularization, which is known to encourage sparsity in the coefficients of regression models \citep{Tibshirani1996}, they propose a ReLU stability regularization, which aims at increasing the number of neurons that can be determined active or inactive a priori. Thus, the number of binary variables necessary to model the ANN is reduced, which leads to smaller optimization problems, and in turn to a significant speedup in the verification problems.
The main disadvantage of methods from this first category is that in some applications, it may not always be feasible to train a dedicated neural network surrogate specifically for optimization. In this case, optimization algorithms have to work with networks that are trained, for instance, with simulation in mind, and it is not possible to specify desirable network dimensions and training methods.

The second category therefore includes methods that a) modify existing ANNs after the training phase to improve their properties and b) obtain tighter bounds in existing formulations for ReLU ANNs. 
Among others, this category include different compression methods (e.g., weight pruning) and optimization-based bound tightening (OBBT). Pruning is usually done via the removal of connections of neurons that have small weights \citep{Cacciola2024} and results in smaller networks with approximately the same functional relationship. Several papers found that models may be compressed without significant loss in accuracy \citep{Han2015,Suzuki2020a}. Exact compression methods, i.e., methods that keep the functional relationship described by the ANN intact, are described in \citet{Kumar2019,Serra2020} for neural networks with ReLU activation. By investigating the bounds of each neuron, smaller networks can be obtained if it can be determined a priori that the input to a neuron is non-negative or non-positive for inputs in the relevant input domain. In this case, no variables have to be added to model the maximum operator in the activation function. If such determinations can be made for all neurons in a layer, then the whole layer can be removed and merged with the subsequent layer via matrix multiplication and addition of the biases. This results in optimization problems with fewer optimization variables and hence in a computational speedup. More recently, theoretical parallels with tropical geometry have been used to simplify neural networks with ReLU activation \citep{Smyrnis2020}. OBBT procedures form a second pillar of this category. In \citet{Grimstad2019}, various bound tightening methods for ReLU ANNs and their effect on optimization times are investigated. \citet{Badilla2023} considers LP-based and MILP-based bound tightening for ReLU ANNs. They analyse the trade-offs of computing tighter bounds via a more expensive bound tightening method and the benefits of the tighter bounds in the subsequent optimization problem.

In addition to the aforementioned methods to speed up the solution process of general MILP solvers, the development of solution heuristics is an active field of research. \citet{Tong2024} propose a heuristic which performs a local search by traversing neighboring linear regions and solving LP subproblems in each linear region. The authors show that for the case of minimizing the output of an ANN and finding adversary examples, this heuristic outperforms general MILP solvers, e.g., Gurobi, especially for deeper networks.

Several software packages have been published that facilitate the embedding of neural networks and other machine learning models into larger optimization problems. \texttt{OMLT}~\citep{Ceccon2022} is a Python package which supports neural networks and gradient boosted trees, and sets up the variables and constraints in the optimization environment \texttt{Pyomo} \citep{Bynum2021,Hart2011}. \texttt{Pyomo} offers interfaces to different solvers, e.g., Gurobi \citep{gurobi}, with which the problems can be solved. Alternatively, \texttt{gurobi-machinelearning} or \texttt{PySCIPOpt-ML} \citep{Turner2024} can be used to translate trained regression models, including neural networks to MIP formulations and solve them with \texttt{Gurobi} and \texttt{SCIP}, respectively. All of these options support models trained by different machine learning backends, including \texttt{Keras} \citep{Chollet2015} or \texttt{PyTorch} \citep{PyTorch}. A further alternative is the software package \texttt{reluMIP} \citep{reluMIP2021}, which has interfaces to both \texttt{Pyomo} and \texttt{Gurobi}, but only supports models trained using \texttt{Keras} or \texttt{TensorFlow}.

\medskip
\noindent \textbf{Contributions and Outline.}
We provide a survey of popular and novel approaches to improving the computational efficiency of optimization with embedded feed-forward ANNs with ReLU activation functions.
In addition and for the first time to our knowledge, we quantify the impact of these methods systematically.
We evaluate and compare all of them on the same benchmark problems, using \texttt{Gurobi}, an analysis of big-M coefficients, and a novel method to calculate the number of piecewise-linear regions.
Although the effects of regularization and dropout on the training of ANN and redundancy of the obtained mathematical model have received much attention in the literature, we study this effect systematically in the context of embedded optimization.

The rest of the paper is structured as follows. In \cref{sec:methods}, we state several methods that have been proposed in the literature to facilitate optimization with embedded neural networks. In addition, we introduce an equivalent transformation of ReLU neural networks that reduces the magnitude of big-M coefficients in their MILP formulation. In \cref{sec:results} we evaluate the influence of these methods on the performance of optimization algorithms in numerical studies. We conclude with a discussion of the findings and possible future lines of research in \cref{sec:discussion}.

\section{Methods}\label{sec:methods}

In this section we discuss several methods that have an impact on the overall performance of a MINLP solver such as \texttt{Gurobi}, when applied to optimization problems of type \eqref{prob:embedded}.
We start by introducing two measures of complexity in this context in Section~\ref{subsec:measures}, namely the number of regions partitioning the input domain in which the function $h(x)$ has identical linear output behavior, and the number of stable ReLU neurons.

Then we examine methods that are applicable to trained ANNs. In Section~\ref{subsec:boundtightening} bound tightening approaches for the optimization problem are presented. In Section~\ref{sec:scaling} we propose a novel scaling method that improves the $\ell^1$ regularization term of a pre-trained network without changing its encoded function. This method can be used after completed training of the ANN (a posteriori) and before the optimization is started (a priori).

In Sections~\ref{sec:regularization}, \ref{sec:clipped}, and \ref{sec:dropout} we investigate modifications to the training of the ANN, in particular regularization of training weights, clipped ReLU formulations, and the use of dropout during training.


\subsection{Measures of complexity of ReLU ANNs}\label{subsec:measures}

While the solution of mixed-integer optimization problems is difficult ($\mathcal{NP}$-complete) in general, it is well known that the number of optimization variables and the tightness of relaxations of the integer variables have a major impact on computational runtimes. In the context of embedded ANNs, we shall consider two particular indicators of complexity.

\subsubsection{Number of linear regions of ReLU networks}

ReLU ANN describe piecewise affine-linear functions \citep{Grigsby2022}. Therefore, the network partitions the input domain $\mathcal{X} \subseteq \mathbb{R}^{n_x}$ into regions in which $h(x)$ is affine linear. These regions are typically called \textit{linear regions}. The bounds on the number of linear regions of a neural network with given depth and width was investigated in \citet{Montufar2014} and later improved on in \citet{Raghu2017}. In general, the number of linear regions of a neural network corresponds to the number of feasible activation patterns in \eqref{eq:bigM}, i.e., the binary decisions whether a neuron is on or off for all neurons in the neural network. Thus, it is an important statistic when considering the complexity of optimizing over neural networks, e.g., in branch-and-bound frameworks, where the variables to branch on represent active or inactive neurons. 

\subsubsection{Number of stable ReLU neurons}

The number of variables a branch-and-bound method has to branch on is an important statistic for estimating the complexity of the optimization problem. In ReLU networks, the variables to branch on are the binary variables $z$ in \eqref{eq:bigM} of every neuron in the network. However, if a neuron can be identified as stable, no binary variable has to be added to model the neuron. To identify stable neurons, their pre-activation bounds are used. The neuron $i$ in layer $j$ is called stably active if $L^{(j)}_i > 0$ and stably inactive if $U^{(j)}_i < 0$, for $j \in [J]$ and $i \in [n_j]$ for all inputs in the input domain $\mathcal{X} \subseteq \mathbb{R}^{n_x}$. 

A regularization to induce ReLU stability was proposed in \citet{Xiao2019} to speed up verification of ReLU networks, whereas in \citet{Serra2020} stable neurons are used to compress neural networks. 
To enumerate the linear regions of a ReLU ANN, we exploit the fact that, within a given linear region, the input and output of each neuron is an affine linear functional in the ANN's overall input space. We use forward sensitivity propagation to calculate the gradient of each neuron's regional input functional and simultaneously perform a forward evaluation of the linearized ANN at the input space's coordinate origin to determine each affine input functional's output shift. With both gradient and shift, we can determine a hyperplane in input space along which the neuron's ReLU activation would switch. We then construct a linear equation system that describes the intersection of halfspaces within which all neurons would retain their current activation pattern. We add the bounds of the input domain to this equation system to ensure boundedness of the linear region. We then use a variant of the \texttt{QuickHull} algorithm~\citep{Barber1996} via the \texttt{SciPy} library~\citep{SciPy2020} to reduce this equation system into an irredundant one and to determine the vertices of the linear region. This also reveals information on which neurons define the facets of the linear region, which means that we can jump across these facets to adjacent regions by switching the activity of those neuron's activation functions. Assuming that there is no facet along which two neurons switch simultaneously, this allows us to enumerate all linear regions that intersect the input domain. We can detect the edge case of two neurons switching simultaneously because it would cause us to enter a region with an empty interior. We do not observe this behavior. 

\subsection{bound tightening} \label{subsec:boundtightening}

Calibrating the big-M coefficients in MILP formulations is crucial for performance of optimization algorithms. bound tightening plays an important role in this context. With ReLU ANNs, there are different ways to compute the big-M coefficients.
\subsubsection{Interval arithmetic}\label{subsec:ia_bounds}
In the presence of input bounds $L^{(0)} \leq x^{(0)} \leq U^{(0)}$ with $L^{(0)}, U^{(0)} \in \mathbb{R}^{n_x}$, big-M coefficients of formulation \eqref{eq:bigM} can be computed via interval arithmetic (IA).
\begin{align}
    L^{(k)}_i &= \sum_{j=1}^{n_{k-1}} \min \{ W_{i,j}^{(k)} L_j^{(k-1)}, W_{i,j}^{(k)} U_j^{(k-1)} \} + b^{(k)}_i, \quad k \in [J],  i \in [n_k] \\
    U^{(k)}_i &= \sum_{j=1}^{n_{k-1}} \max \{ W_{i,j}^{(k)} L_j^{(k-1)}, W_{i,j}^{(k)} U_j^{(k-1)} \} + b^{(k)}_i, \quad k \in [J],  i \in [n_k] 
\end{align}
This forward propagation yields valid bounds. However, it ignores the fact that the activation of neurons, i.e., whether they are on the left or right arm of the ReLU function, is not independent between neurons. This results in overly relaxed approximations of the actual bounds. As a result, there is typically an exponential increase of big-M coefficients with increasing depth. This behaviour is exemplified in \vref{fig:IA_bounds}.

\subsubsection{LP-based bound tightening}\label{subsec:lr_bounds}

The bounds from \cref{subsec:ia_bounds} can be tightened by taking advantage of dependencies between the neurons as well as potentially existing bounds on the output of the neural network $L^{(J)} \leq x^{(J)} \leq U^{(J)}$. This is achieved by solving two auxiliary optimization problems per neuron, minimizing and maximizing, respectively, the pre-activation value of each neuron. The optimization problem for computing tighter bounds for neuron $k$ in layer $j$, with $j \in [J],\,k \in [n_k]$ in its general form as an MILP reads

\begin{equation}\label{prob:obbt}
    \begin{aligned}
    \min_{x,z}\ & W^{(j)}_k x^{(j-1)} + b^{(j)}_k \\ 
    \textrm{s.t.}\ & \begin{alignedat}[t]{3}
            x^{(j)}_i & \geq 0, \quad && j \in [J],\, i \in [n_j] \\
            x^{(j)}_i & \geq W^{(j)}_i x^{(j-1)} + b^{(j)}_i, \quad && j \in [J],\, i \in [n_j] \\
            x^{(j)}_i & \leq W^{(j)}_i x^{(j-1)} + b^{(j)}_i - L^{(j)}_i (1-z^{(j)}_i), \quad && j \in [J],\, i \in [n_j] \\
            x^{(j)}_i & \leq U^{(j)}_i z^{(j)}_i, \quad && j \in [J],\, i \in [n_j] \\
            x^{(0)}_i  & \leq U^{(0)}_i,        \quad && i \in [n_x]\\
            x^{(0)}_i  & \geq L^{(0)}_i,        \quad && i \in [n_x]\\
            z^{(j)}_i & \in \{0,1\}. \quad && j \in [J],\, i \in [n_j]
        \end{alignedat}
    \end{aligned}
\end{equation}
Solving \eqref{prob:obbt} yields a valid lower bound $L^{(j)}_k$, while the corresponding upper bound $U^{(j)}_k$ is computed by maximizing instead of minimizing in \eqref{prob:obbt}. In order to reduce the computational effort, typically the LP relaxation of formulation \eqref{prob:obbt} is considered. Hence, the auxiliary problems are linear programs (LPs) and can be solved efficiently. Solving the MILP directly is considered in \citet{Badilla2023,Grimstad2019}. However, the reduction in computational effort in subsequent optimization is quickly outweighed by the effort spent on solving the bound tightening MILPs. Therefore, we only consider the LP-based bound tightening procedure in this paper. One degree of freedom when performing bound tightening is the ordering of variables for which bounds are tightened. As the direction of bound propagation is from the input to the output layer, this is also the natural order to perform the tightening. However, within each layer the order may be chosen arbitrarily. Different methods to choose this order are discussed in \citet{Rossig2021}. However, they do not find any advantage of more advanced methods over a simple, fixed ordering of variables. Therefore, in this contribution, we apply bound tightening in a fixed ordering of variables.

\subsection{A posteriori scaling of ReLU ANNs} \label{sec:scaling}
Weights of neural networks are not uniquely determined by the training process and the training data, i.e., there are different realizations of weights and biases that define the same functional relationship of input and output. This observation can be exploited to design algorithms that transform a trained neural network into a functionally identical network with some desired property. This could be, e.g., a lower norm of the weight matrices. With the input bounds remaining unchanged, this would lead to a reduction of big-M coefficients, which could be beneficial in subsequent optimization problems.

\begin{figure}
    \centering
    \includegraphics[width=.6\linewidth]{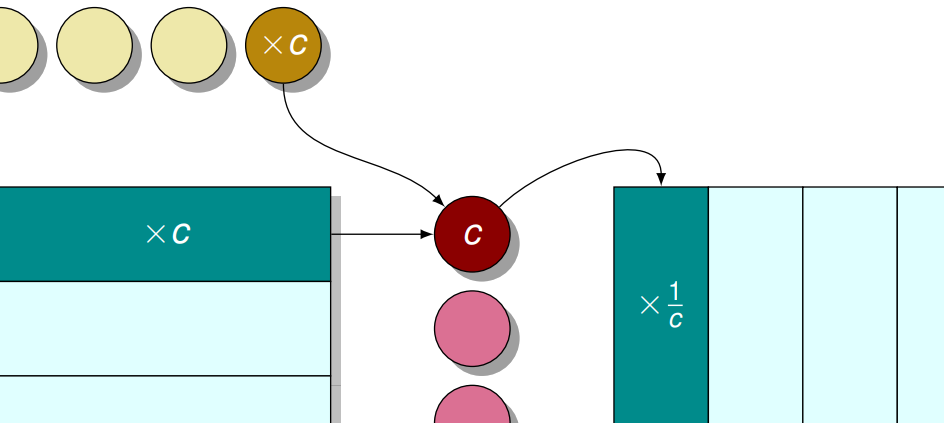}
    \caption{Equivalent scaling of ReLU ANNs. Scalar factor $c$ is multiplied row-wise to weight matrix and corresponding bias of current layer, resulting in a scaling of the output of the neuron by a factor of $c$. To compensate this, the weight matrix in the subsequent layer needs to be multiplied column-wise with the reciprocal of $c$.}
    \label{fig:rescale}
\end{figure}

In case of the ReLU activation function, one can exploit its positive homogeneity. For a single neuron $i$ in layer $k$, with $k \in [J], \ i \in [n_k]$ and a scalar $c^{(k)}_i > 0$ it holds, that
\begin{align}
    \myReLU{c^{(k)}_i \left( W^{(k)}_i x^{(k-1)} + b_i \right)} = c^{(k)}_i \cdot \myReLU{ W^{(k)}_i x^{(k-1)} + b_i},
\end{align}
with $W^{(k)}_i \in \mathbb{R}^{1 \times n_{k-1}}$ being the $i$-th row of the weight matrix in layer $k$. 
To ensure the functional equivalence of the neural network, the $i$-th column of the weight matrix of layer $k+1$, corresponding to the scaled neuron $i$ in layer $k$, needs to be multiplied with the reciprocal of $c^{(k)}_i$. As all neurons of the neural network may be scaled, all weight matrices except the first and the last are scaled with the ratio of the two scaling factors of their surrounding layers. As the bias is not multiplied with the output from the previous layer, no multiplication with the reciprocal is needed. In the final layer $J$, no more new scaling factors may be introduced as they can no longer be compensated in subsequent layers. Therefore, only the scaling of layer $J-1$ is compensated by multiplying $W^{(J)}$ with the reciprocals of the scaling factors of the penultimate layer. For any set of scaling factors $c^{(k)}_i > 0,\,k \in [J], i \in [n_k]$, scaled weights and biases  $\tilde{W}$ and $\tilde b$, computed as 
\begin{equation}
    \begin{alignedat}{3}
        \tilde{W}_{i,j}^{(1)} &= W_{i,j}^{(1)} \cdot c_i^{(1)}, \quad && i \in [n_1],\, j \in [n_x] \\       
        \tilde{W}_{i,j}^{(k)} &= W_{i,j}^{(k)} \cdot \frac{c_i^{(k)}}{c_j^{(k-1)}}, \quad && k \in \{2, \ldots, J-1\},\, i \in [n_k],\, j \in [n_{k-1}],\\
        \tilde{W}_{i,j}^{(J)} &= W_{i,j}^{(J)} \cdot \frac{1}{c_j^{(J-1)}}, \quad && i \in [n_J],\, j \in [n_{J-1}] \\
        \tilde{b}_i^{(k)} &= b_i^{(k)} \cdot c_i^{(k)}, \quad && k \in [J-1], \, i \in [n_k]
    \end{alignedat}
\end{equation}
define functionally equivalent neural networks.
This basic idea of an equivalent transformation of ReLU networks via scaling one layer and compensating the effect of scaling in the next layer is illustrated in \vref{fig:rescale}. 

The scaling factors $c_i^{(k)}$ can be chosen arbitrarily. However, we can specifically choose them such that the resulting network has favorable properties. We propose formulating an optimization problem to obtain scaling factors that minimize the absolute value of the scaled weights $\tilde{W}$ and biases $\tilde{b}$. This lower norm of the weights then yields lower big-M coefficients. As noted in \cref{subsec:ia_bounds}, these are determined solely\footnote{While bounds on the output of the network can be propagated backwards through the network and thus influence the big-M coefficients \citep{Grimstad2019}, we refer only to the big-M coefficients derived via interval arithmetic and forward propagation of input bounds as explained in \cref{subsec:ia_bounds}.} by the input bounds and the magnitude of weights and biases. Hence, this approach can be applied to networks that were not initially trained with regularization in order to generate an equivalent neural network with lower big-M coefficients. Of course, other effects of regularization, e.g., weight sparsity, cannot be obtained by this method. The proposed optimization problem is
\begin{equation}
    \begin{aligned}
        \min_{c}\ & \sum_{i=1}^{n_1} \sum_{j=1}^{n_x} \lvert W_{i,j}^{(1)} \rvert \cdot c_i^{(1)} +  \sum_{k=2}^{J-1} \sum_{i=1}^{n_k} \sum_{j=1}^{n_{k-1}}  \lvert W_{i,j}^{(k)} \rvert \cdot \frac{c_i^{(k)}}{c_j^{(k-1)}} \\*
        & + \sum_{k=1}^{J-1} \sum_{i=1}^{n_k} \lvert b_i^{(k)} \rvert \cdot c_i^{(k)} + \sum_{i=1}^{n_J} \sum_{j=1}^{n_{J-1}} \lvert W_{i,j}^{(J)} \rvert \cdot \frac{1}{c_i^{(J)}} \\
        \textrm{s.t.}\ & \begin{alignedat}[t]{3}
            c_i^{(k)} & > 0, \quad && k \in [J],\, i \in [n_k] \\
            c & \in \bigtimes_{k=1}^J \mathbb{R}^{n_k} \quad
        \end{alignedat}
    \end{aligned}
\end{equation}
This optimization problem is not trivial to solve directly because it involves fractions and strict inequality constraints. However, because all $c_i^{(k)}$ have to be strictly positive, we can convert it into a convex optimization problem on a closed set by replacing each $c_i^{(k)}$ with its logarithm. Each summand in the objective function then becomes an evaluation of the exponential function, multiplication becomes addition, and division becomes subtraction. With the logarithm of $c^{(k)}_i$  referred to as $\tilde{c}^{(k)}_i$, the transformed optimization problem reads
\begin{equation}\label{prob:scaling}
    \begin{aligned}
        \min_{\tilde{c}}\ & \sum_{i=1}^{n_1} \sum_{j=1}^{n_x} \exp \left( \log \left( \lvert W_{i,j}^{(1)} \rvert\right) + \tilde{c}_i^{(1)} \right)   + \sum_{k=2}^{J} \sum_{i=1}^{n_k} \sum_{j=1}^{n_{k-1}} \exp \left( \log \left( \lvert W_{i,j}^{(k)} \rvert \right)+ \tilde{c}_i^{(k)} - \tilde{c}_j^{(k-1)} \right)\\*
        & + \displaystyle \sum_{k=1}^{J} \sum_{i=1}^{n_k} \exp \left( \log \left( \lvert b_i^{(k)} \rvert \right) +  \tilde{c}_i^{(k)} \right) + \sum_{i=1}^{n_J} \sum_{j=1}^{n_{J-1}}  \exp \left( \log\left(\lvert W_{i,j}^{(J)} \rvert \right) - \tilde{c}_i^{(J)} \right) \\
        \textrm{s.t.}\ & \tilde{c} \in \bigtimes_{k=1}^J \mathbb{R}^{n_k}
    \end{aligned}
\end{equation}

\subsection{Regularization} \label{sec:regularization}

The objective function for training neural networks typically consists of two terms. The first 
accounts for the mismatch between prediction and data, while the second term aims at preventing overfitting and thus allowing for a better generalization of the model to unseen data. With $W \in \mathbb{R}^d$ denoting the vector of all weights and biases and $N \in  \mathbb{N}$ representing the number of training samples of inputs and outputs $(x_i, y_i), \, i \in [N]$, the objective reads
\begin{equation}\label{prob:training}
    \begin{aligned}
        \min_{W}\ & \frac{1}{N} \sum_{i=1}^{N} \left( h(x_i) -y_i\right)^2 + \lambda \Omega(W)
    \end{aligned}
\end{equation}
Popular choices for the regularization term $\Omega : \mathbb{R}^d \mapsto \mathbb{R}$ are the penalization of large magnitudes of weights and biases by using some vector norm, e.g., $\Omega(W) = \| W \|_p$, with typically $p=1$ and $p=2$. Typical ways to measure the generalization performance of a model is to compute the mean absolute percentage error (MAPE) defined as 
\begin{align}
    \text{MAPE}\left(\hat{y}, y\right) = \frac{1}{n} \sum_{i=1}^{n} \frac{|\hat{y}_i - y_i|}{\max\{\varepsilon, |y_i|\}}
\end{align}
for predictions $\hat{y}$ on the test dataset.

While it is known that $\ell^1$ regularization leads to sparser regression models \citep{Tibshirani1996}, \citet{Xiao2019,Serra2020} found that applying $\ell^1$ regularization also increased ReLU stability, i.e., the percentage of stable neurons. The authors of \citet{Xiao2019} also propose a dedicated ReLU stability regularization \eqref{eq:stability_regularization}, which penalizes the sign differences in the pre-activation bounds of each neuron, thus encouraging stability.
\begin{align}\label{eq:stability_regularization}
    \Omega_{\text{RS}}(W) = - \sum_{i=1}^{J} \sum_{j=1}^{n_i} \text{sign}(U_j^{(i)}) \cdot \text{sign}(L_j^{(i)})
\end{align}
For practical purposes, a smooth reformulation of \eqref{eq:stability_regularization} is used, and \citet{Xiao2019} show that verification problems of neural networks trained using this regularization can be solved faster than with $\ell^1$ regularization due to a higher number of stable neurons. In this paper, we will however focus on investigating the effect of varying levels of $\ell^1$ regularization on the performance of optimization algorithms, as it is one of the most commonly used types of regularization.

\subsection{Clipped ReLU} \label{sec:clipped}
One of the reasons why big-M coefficients in ReLU networks increase quickly with increasing network depth is that the ReLU activation function is unbounded. A variation of the ReLU function is the clipped ReLU function proposed in \citet{Hannun2014}. In the clipped ReLU function, the output of the function is bounded by an upper value $M \in \mathbb{R}$, i.e.,
\begin{align}\label{eq:relu_m}
    \myClippedReLU[M]{x} = \max\bigl\{0,\, \min\{M,\, x\}\bigr\}
\end{align}
Using standard disjunctive programming notation, the feasible set of $x_i^{(j)} = \myClippedReLU[M]{ W^{(j)}_i x^{(j-1)} + b_i}$ can be written as
\[
\begin{bmatrix}
    x_i^{(j)} = 0 \\
    W^{(j)}_i x^{(j-1)} + b_i \leq 0
\end{bmatrix}
\vee 
\begin{bmatrix}
    x_i^{(j)} = W^{(j)}_i x^{(j-1)} + b_i\\
    0 < W^{(j)}_i x^{(j-1)} + b_i < M 
\end{bmatrix}
\vee
\begin{bmatrix}
    x_i^{(j)} = M \\
    W^{(j)}_i x^{(j-1)} + b_i \geq M 
\end{bmatrix}
\]
We formulate a big-M relaxation of this feasible set as
\begin{align}
    \begin{split}
        x_i^{(j)} & \geq 0, \\
        x_i^{(j)} & \leq M z_{1_i}^{(j)}, \\
        x_i^{(j)} & \leq U_i^{(j)} z_{1_i}^{(j)}, \\
        x_i^{(j)} & \leq W^{(j)}_i x^{(j-1)} + b_i - L_i^{(j)} \cdot (1-z_{1_i}^{(j)}),\\
        x_i^{(j)} & \geq M z_{2_i}^{(j)}, \\
        x_i^{(j)} & \geq W^{(j)}_i x^{(j-1)} + b_i - (U_i^{(j)}-M) z_{2_i}^{(j)}, \\
        z_{1_i}^{(j)}, z_{2_i}^{(j)} & \in \{0,1\},\\
        z_{1_i}^{(j)} & \geq z_{2_i}^{(j)},
    \end{split}
    \label{eq:bigM_clipped}
\end{align}
similar to the formulation suggested in a preprint version of \citet{Anderson2020}.
This formulation comes at the cost of an additional binary variable compared to the standard big-M formulation \eqref{eq:bigM}. If both binary variables are zero, the neuron is inactive and $x_i^{(j)}=0$. In the case $z_{1_i}^{(j)}=1, z_{2_i}^{(j)}=0$, the neuron is active and $0 \leq x_i^{(j)} =  W^{(j)}_i x^{(j-1)} + b_i \leq M$. If both binary variables are non-zero, the neuron's output is limited by the threshold $M$. 

\subsection{Dropout} \label{sec:dropout}

Dropout is a technique applied during training proposed in \citet{Srivastava2014} to prevent overfitting the data by randomly turning off a percentage of the neurons in some or all layers. Therefore, redundancies have to be established in the neural network to achieve an adequate accuracy. There is empirical evidence that neural networks trained with dropout have more linear regions \citep{Zhang2020a} than those trained without. Hence, in contrast to the aforementioned methods, it is expected that applying dropout during training leads to more complex neural networks which makes optimizing over them more difficult.
We will thus apply dropout as an antithesis to validate our conjecture that the runtime of MINLP solvers increases for more redundant and decreases for less redundant ANN models.

\section{Numerical results}\label{sec:results}

In the \cref{sec:methods}, we have enumerated some methods to formulate, train, and scale feed-forward neural networks with ReLU activation (or variations thereof), as well as to tighten their relaxation prior to optimization through bound tightening. In this section, we evaluate how these methods affect global optimization performance. In order to do so, we train neural networks as surrogates for several non-convex benchmark functions and compare solver performance with various post-processing steps.

We first present numerical results on relevant characteristics of ReLU ANNs in the context of optimization. These include their expressive power as measured by the number of linear regions they define and the percentage of stable neurons that can be determined from the pre-activation bounds, introduced in the beginning of \cref{sec:methods}.  We count only those linear regions that intersect the relevant input domain of each function.

We show how the methods presented in \cref{sec:methods} impact these quantities and improve the performance of optimization algorithms. For this, we restrict ourselves to minimizing the output of feed-forward ReLU ANNs, i.e., the optimization problem we solve reads
\begin{equation}\label{prob:minANN}
    \min_{x} \ h(x)
\end{equation}
where $h \colon \mathbb{R}^{{n_x}} \mapsto \mathbb{R}$ is the trained neural network. 
The benchmark functions we consider for approximation and subsequent minimization are:

\begin{figure}[t]
    \centering
    \begin{subfigure}{.32\linewidth}
        \includegraphics[width=\linewidth]{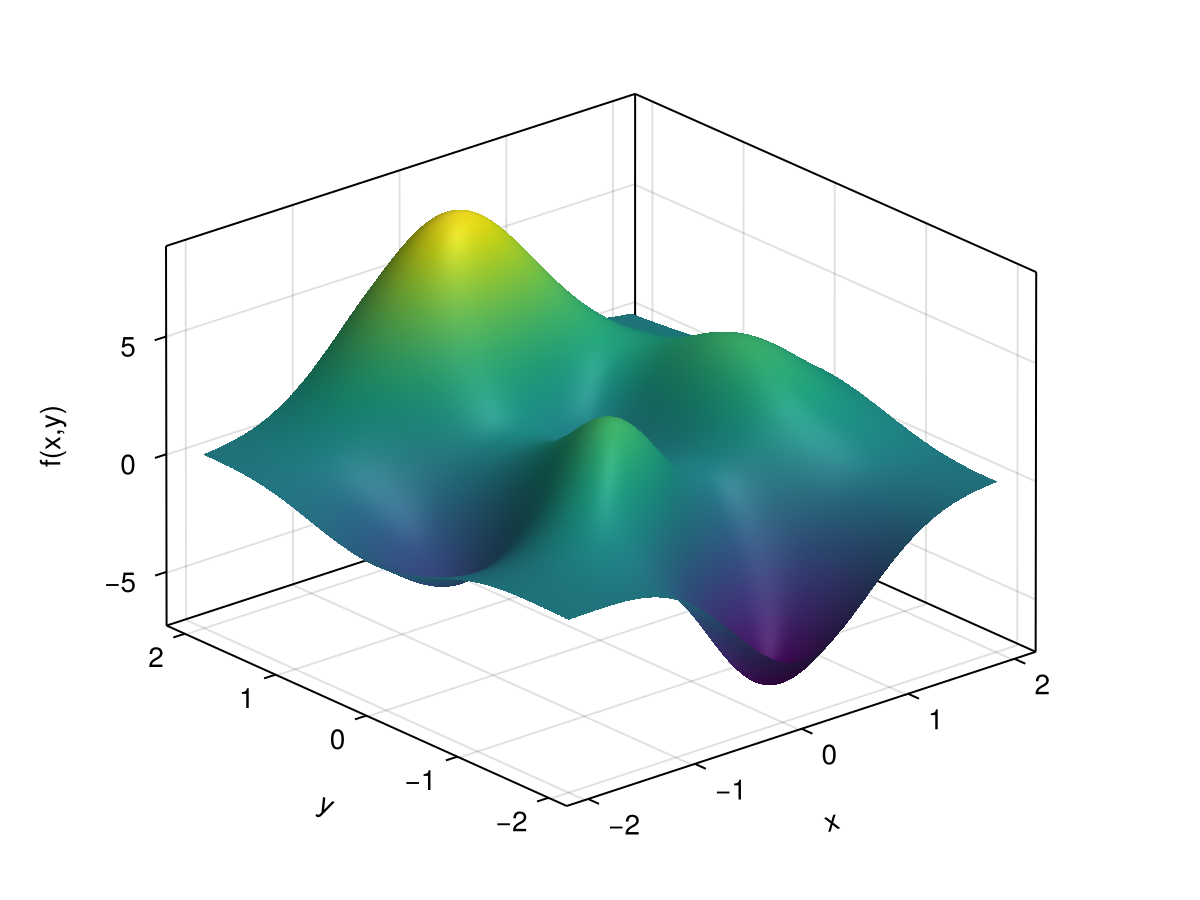}
        \subcaption{Peaks function \eqref{fun:peaks}}\label{fig:peaks-function}
    \end{subfigure}
    \begin{subfigure}{.32\linewidth}
        \includegraphics[width=\linewidth]{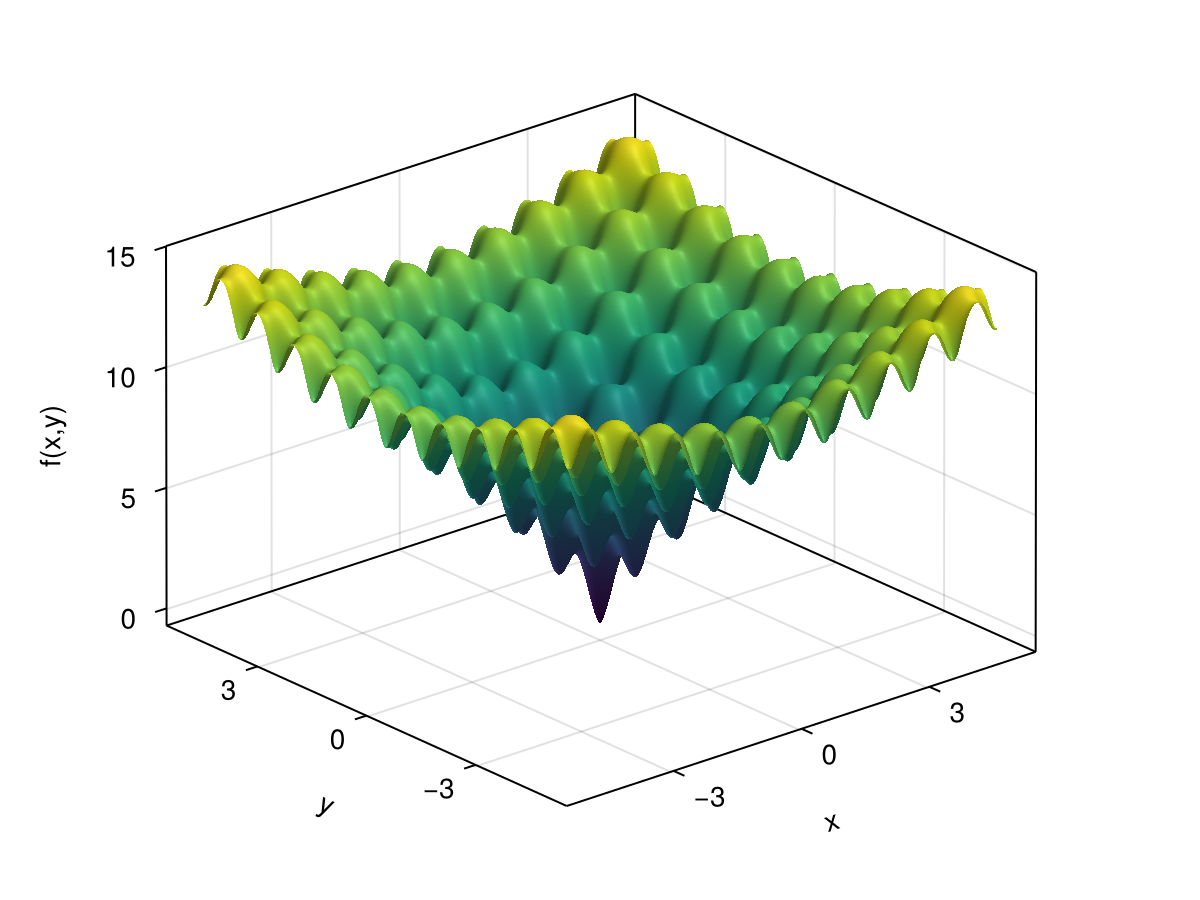}
        \subcaption{Ackley's function \eqref{fun:ack}}\label{fig:ack-function}
    \end{subfigure}
    \begin{subfigure}{.32\linewidth}
        \includegraphics[width=\linewidth]{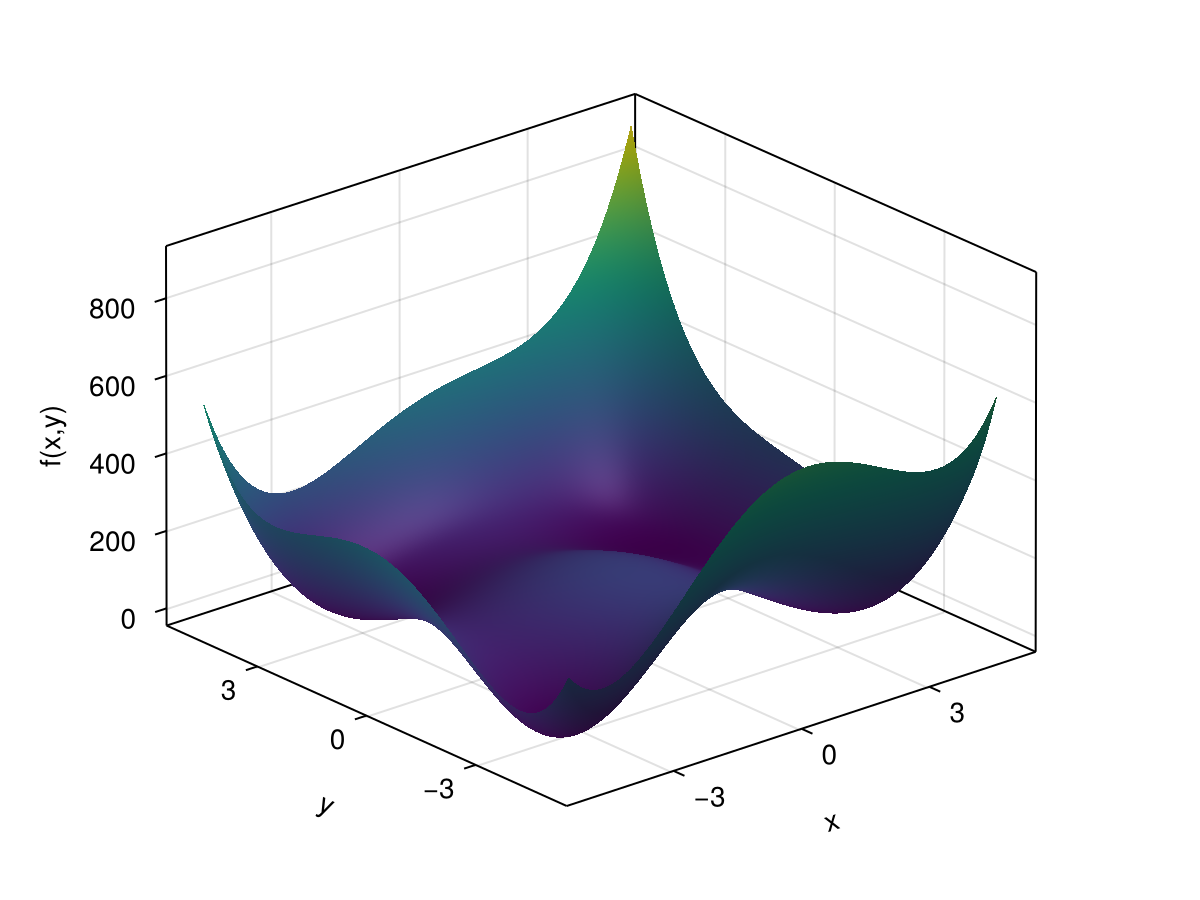}
        \subcaption{Himmelblau's function \eqref{fun:him}}\label{fig:him-function}
    \end{subfigure}
    \caption{Surface plots of the benchmark functions for surrogate model training and optimization.}
    \label{fig:different_functions}
\end{figure}

\begin{enumerate}
    \item The Peaks function $f_{\text{peaks}} \colon \mathbb{R}^2 \mapsto \mathbb{R}$ is given by
    \begin{align}
        \begin{split}
        f_{\text{peaks}}(x,y) &= - 3 (1-x)^2 \exp\bigl(-x^2 - (y+1)^2\bigr)   - 10 \Bigl( \frac{x}{5} - x^3 - y^5\Bigr) \exp(-x^2 - y^2) \\
            & \quad - \frac{1}{3}  \exp\bigl(-(x+1)^2 - y^2\bigr).
        \end{split}
        \label{fun:peaks}
    \end{align}
    It is commonly used as a benchmark function, e.g., in \citet{Schweidtmann2019a} and has multiple local minima and maxima on the domain $x,y \in [-2,2]$. The global minimum is $(0.228, -1.626)$ with objective value $-6.551$. The function is depicted in \vref{fig:peaks-function}.

    \item Ackley's function $f_{\text{ackley}} \colon \mathbb{R}^2 \mapsto \mathbb{R}$ is defined by
    \begin{align}
        \begin{split}
            f_{\text{ackley}}(x,y) & = -20 \cdot \exp \left(-\frac{1}{5}\sqrt{\frac{1}{2}(x^2+y^2)}\right) \\
                            & \quad  - \exp \left(\frac{1}{2} \bigl(\cos(2\pi x) + \cos(2\pi y)\bigr) \right) + \exp(1) + 20    
        \end{split}
        \label{fun:ack}
    \end{align}
    and is often used as a benchmark function for optimization algorithms. For instance, it is used in \citet{Tsay2021}. It is considered on the domain $x,y \in [-3.5,3.5]$. It is non-convex, has several local minima and one global minimum at $x=y=0$ with objective value $0$. A surface plot is depicted in \vref{fig:ack-function}.

    \item Himmelblau's function $f_{\text{himmelblau}} \colon \mathbb{R}^2 \mapsto \mathbb{R}$ with
    \begin{align}\label{fun:him}
        f_{\text{himmelblau}}(x,y) =  (x^2 + y - 11)^2 + (x + y^2 - 7)^2
    \end{align}
    is considered on the domain $x,y \in [-5,5]$, where it has four equivalent local (and global) minima: $(3.0,2.0)$, $(-2.805, 3.131)$, $(-3.779, -3.283)$ and $(3.584,-1.848)$. All have objective function value $0$. A surface plot is depicted in \vref{fig:him-function}.
\end{enumerate}

In our numerical study, we consider a total of 1080 different neural networks. This number of instances stems from considering the three benchmark functions used for the approximation of \cref{fun:peaks,fun:ack,fun:him} and the specific choices for the hyperparameters of the trained neural networks. These differ both in their width and depth, as well as the activation function and the level of $\ell^1$ regularization applied during training. The specific options for these hyperparameters are given in \Cref{table:hyperparameter}. Using Latin Hypercube sampling, we generated training data of 100,000 samples for the Peaks function \eqref{fun:peaks} and Himmelblau's function \eqref{fun:him}, and 150,000 samples for Ackley's function \eqref{fun:ack}, to account for its higher nonconvexity. For training, we first normalize both input and output data, and reserve 30\% of the data as a test set to evaluate the generalization of the networks. All networks are then trained for 300 epochs using the Adam algorithm \citep{Kingma2017}. To study the effect of scaling and bound tightening on each of the trained networks, we solve problems \eqref{prob:scaling} and \eqref{prob:obbt}, where applicable. As the scaling method is not designed for the clipped ReLU, we can only solve \eqref{prob:scaling} for the 360 instances with standard ReLU activation. We use \texttt{OMLT} \citep{Ceccon2022} to set up the constraints for the ReLU ANNs via \texttt{Pyomo}~\citep{Bynum2021,Hart2011}, and \texttt{Gurobi}~\citep{gurobi}~{v11.0.1} with default options and a time limit of 300 seconds to solve the resulting optimization problems. 
\begin{figure}
    \centering
    \begin{subfigure}{0.49\textwidth}
        \includegraphics[width=\textwidth]{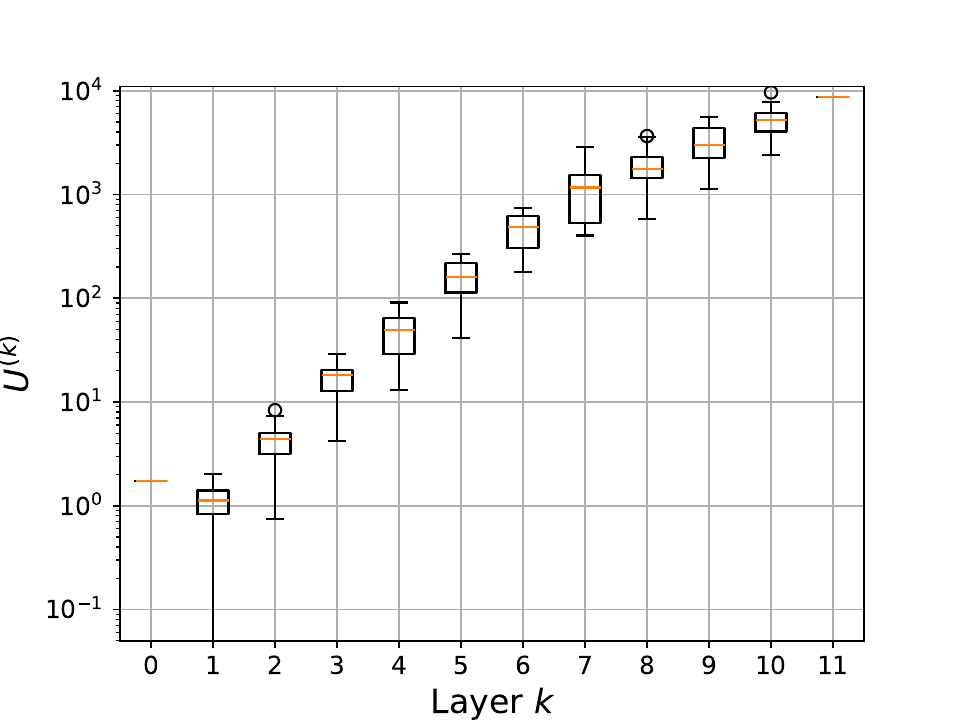}
        \caption{Big-M coefficients $U^{(k)}$ determined via IA for standard ReLU ANN.}
        \label{fig:IA_bounds}
    \end{subfigure}
    \hfill
    \begin{subfigure}{0.49\textwidth}
        \includegraphics[width=\textwidth]{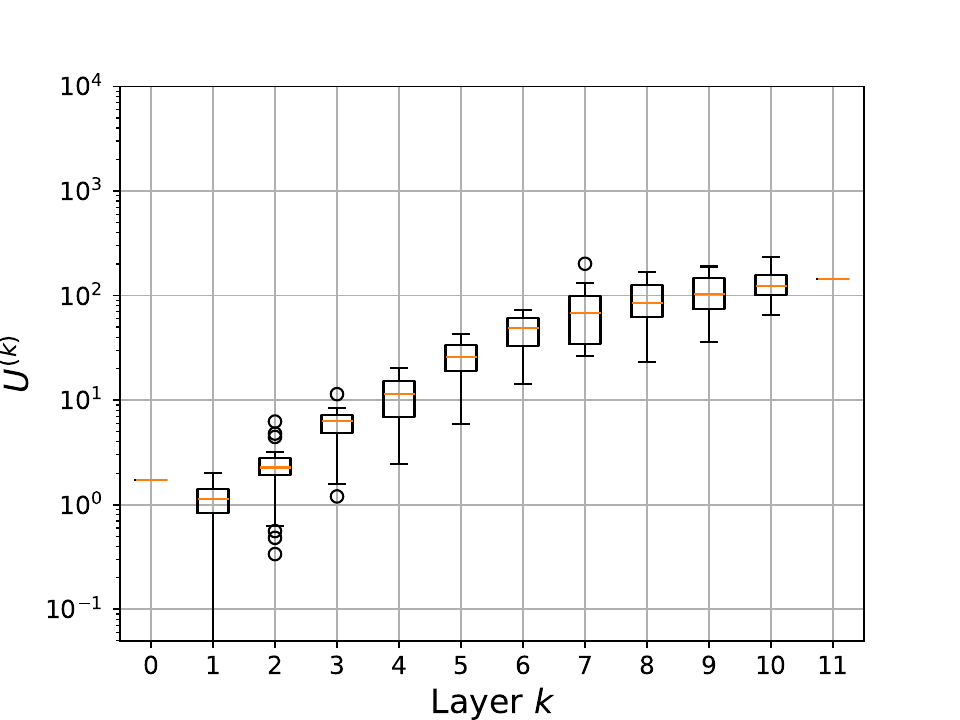}
        \caption{Big-M coefficients $U^{(k)}$ determined via OBBT for standard ReLU ANN.}
        \label{fig:LR_bounds}
    \end{subfigure}
    
    \begin{subfigure}{0.49\textwidth}
        \includegraphics[width=\textwidth]{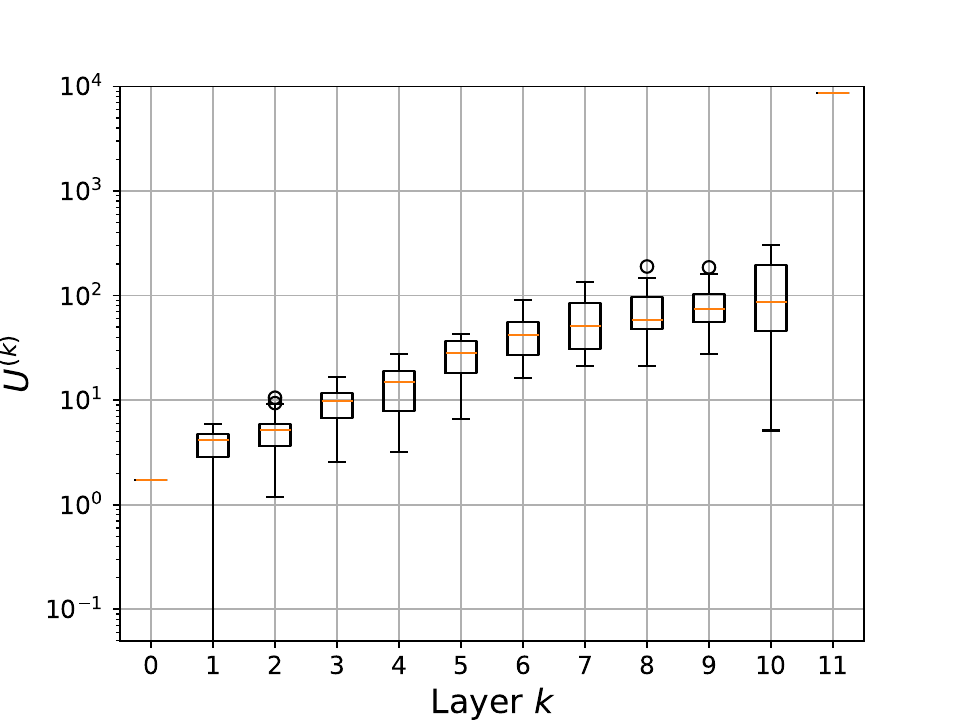}
        \caption{Big-M coefficients $U^{(k)}$ determined via IA for ReLU ANN after ReLU scaling.}
        \label{fig:scaled_IA_bounds}
    \end{subfigure}
    \hfill
    \begin{subfigure}{0.49\textwidth}
        \includegraphics[width=\textwidth]{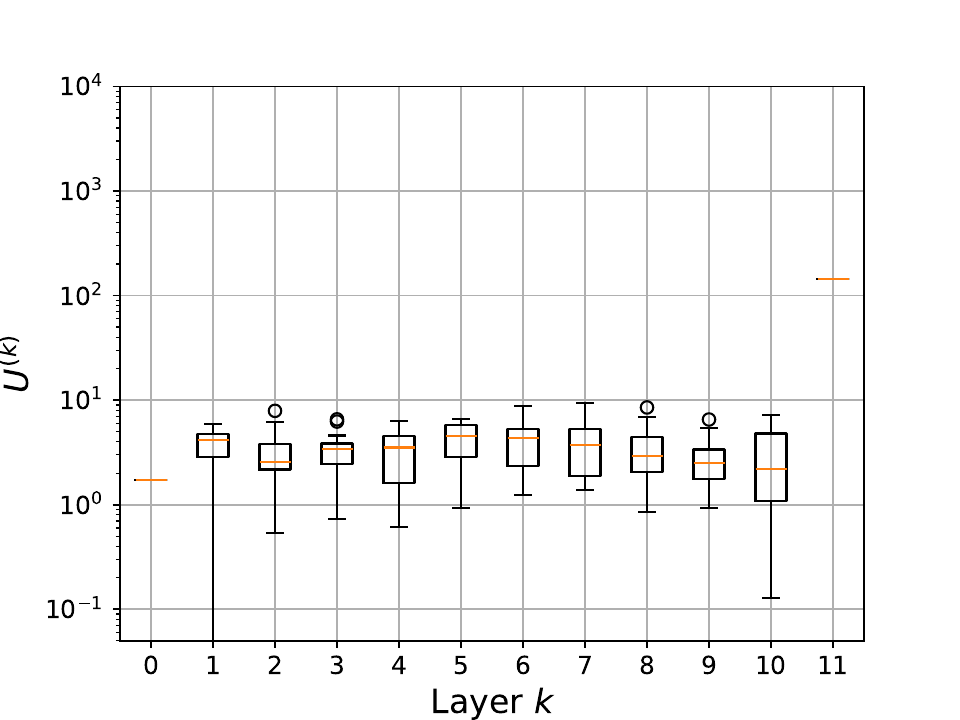}
        \caption{Big-M coefficients $U^{(k)}$ determined via OBBT for ReLU ANN after ReLU scaling.}
        \label{fig:scaled_LR_bounds}
    \end{subfigure}
            
    \caption{Comparison of pre-activation bounds  $U^{(k)}$ for functionally equivalent neural networks with ten hidden layers. The original bounds derived via interval arithmetic shown in \ref{fig:IA_bounds} are characterized by the typical exponential increase due to forward propagation of the input bounds. Solving auxiliary LPs yields tighter bounds, although the exponential increase is still present, as shown in \ref{fig:LR_bounds}. Comparable bounds can be computed via solving the scaling problem \eqref{prob:scaling}, with the distinction that the bounds on the output of the network are equivalent to those derived from interval arithmetic. For the scaled neural network, solving the bound tightening problem \eqref{prob:obbt} in addition yields even tighter bounds on the big-M coefficients in the hidden layers with ReLU activation, as can be seen in \ref{fig:scaled_LR_bounds}, while the output bounds are equivalent to those in \ref{fig:LR_bounds}.}
    \label{fig:big_Ms}
\end{figure}

\begin{table}[ht]
    \centering
    \caption{Hyperparameter options for training of neural networks. Besides varying the depth and width of the networks, we investigate two variants of the clipped ReLU activation \eqref{eq:relu_m} and five levels of $\ell^1$ regularization. All hidden layers have the same dimension.}
    \label{table:hyperparameter}
    \begin{tabular}{l c } 
        \toprule 
        Hyperparameter & Options \\ \midrule
        Hidden Layers & $1,\ldots,10$  \\
        Layer Width &  $25,50$ \\
        Activation & ReLU, ReLU$_2$, ReLU$_5$ \\ 
        $\lambda$ &  $0.0, 10^{-7}, 10^{-6},10^{-5}, 10^{-4}, 10^{-3}$\\ \bottomrule 
    \end{tabular}
\end{table}

\subsection{Effect of OBBT}

\begin{figure}
    \centering
    \begin{subfigure}{.47\linewidth}
        \centering
        Percentage of stable neurons
        \includegraphics[width=\linewidth]{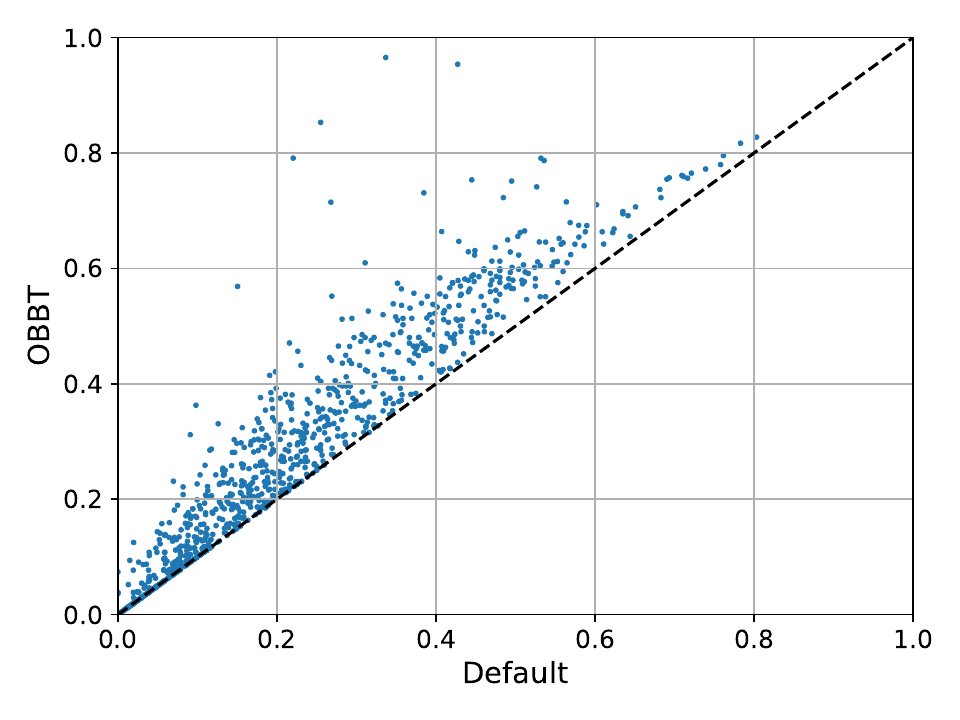}
        \subcaption{Parity plot for percentage of stable neurons compared for bounds from IA and LP-based OBBT.}
        \label{subfig:percentage_fixed_IA_vs_OBBT}
    \end{subfigure}
    \begin{subfigure}{.47\linewidth}
        \centering
        Computational time
        \includegraphics[width=\linewidth]{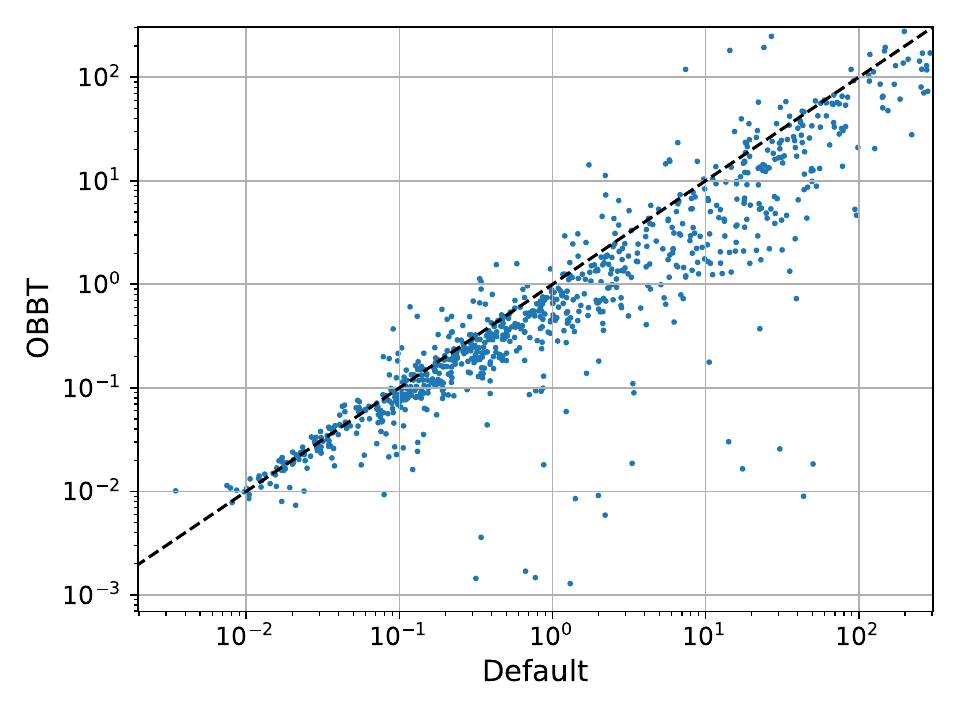}
        \subcaption{Parity plot for computational time compared for bounds from IA and LP-based OBBT.}
        \label{subfig:time_IA_vs_OBBT}
    \end{subfigure}

    \caption{Parity plots comparing percentage of stable neurons and computational times of optimally solved instances of \eqref{prob:minANN} for bounds derived from IA and LP-based OBBT. Solving \eqref{prob:obbt} leads to an increase of 5.5 percentage points in stable neurons on average. This carries over to a reduction in computational time shown in \subref{subfig:time_IA_vs_OBBT}. The ratios of times with OBBT and IA bounds have a geometric mean of 0.57.}
    \label{fig:improvement_obbt}
\end{figure}

For the 1080 trained neural networks, we solve the LP-relaxation of \eqref{prob:obbt} to compute tighter big-M coefficients for formulation \eqref{eq:bigM}, and use them in the optimization problem \eqref{prob:minANN}. The effect on a network with ten hidden layers is illustrated in \vref{fig:big_Ms}. Compared to the IA bounds, there is a reduction in big-M coefficients of the last layer by roughly two orders of magnitude. As \vref{tab:results} shows, OBBT is effective for all trained networks. We assess the reduction in big-M coefficients across all networks by comparing the averaged distances between upper and lower bound $U^{(j)}_k -L^{(j)}_k$ for bounds based on OBBT and IA. Then, over all networks, we calculate the geometric mean over the ratio of these averages. The resulting geometric mean of 0.54 suggests, that, as a rough estimate, OBBT is reducing the big-M coefficients by half. 
As a side effect of these tighter bounds, the percentage of stable neurons increases by 5.5 percent on average. We assess the resulting improvement in computational times by calculating the ratios of the measured  computational times with tightened bounds and those with the original bounds, restricted to instances that were solved to global optimality in both cases. Over these ratios, we again form the geometric mean. With a geometric mean of 0.57, bound tightening brings a significant computational speedup, though it does not substantially increase the number of instances that are solved within the time limit. \vref{fig:improvement_obbt} illustrates the parities of percentage of stable neurons and computational time.

\subsection{Effect of ReLU scaling}

\Vref{fig:big_Ms} illustrates the effects of solving the scaling problem \eqref{prob:scaling} on the big-M coefficients of a neural network with ten hidden layers. The first observation is that the output bounds remain unchanged compared to the original neural network, which is expected as the functional relationship is equivalent. However, the lower $\ell^1$ norm of the weights leads to a reduction in the big-M coefficients for the hidden layers. They are roughly on the same order of magnitude as those obtained via LP-based bound tightening. When both scaling and bound tightening are applied sequentially, the bounds for the hidden layers are tighter than those achieved by OBBT on its own. Also, with the sequential application of scaling and tightening we do not observe any clear sign of an exponential increase in bounds with increasing depth.

\begin{figure}
    \centering
    \begin{subfigure}{.47\linewidth}
        \centering
        Computational time
        \includegraphics[width=\linewidth]{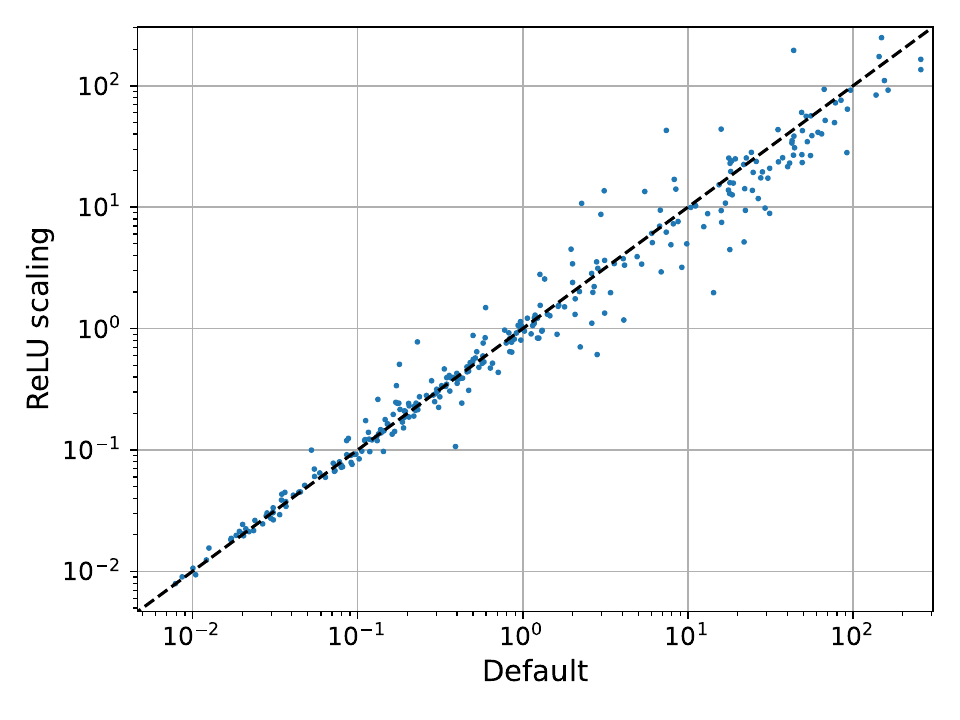}
        \subcaption{Runtime comparison between IA bounds for the baseline network (\enquote{Default}) and IA for the scaled ANN (\enquote{ReLU scaling}).}
        \label{subfig:IA_vs_scaler}
    \end{subfigure}
    \begin{subfigure}{.47\linewidth}
        \centering
        Computational time
        \includegraphics[width=\linewidth]{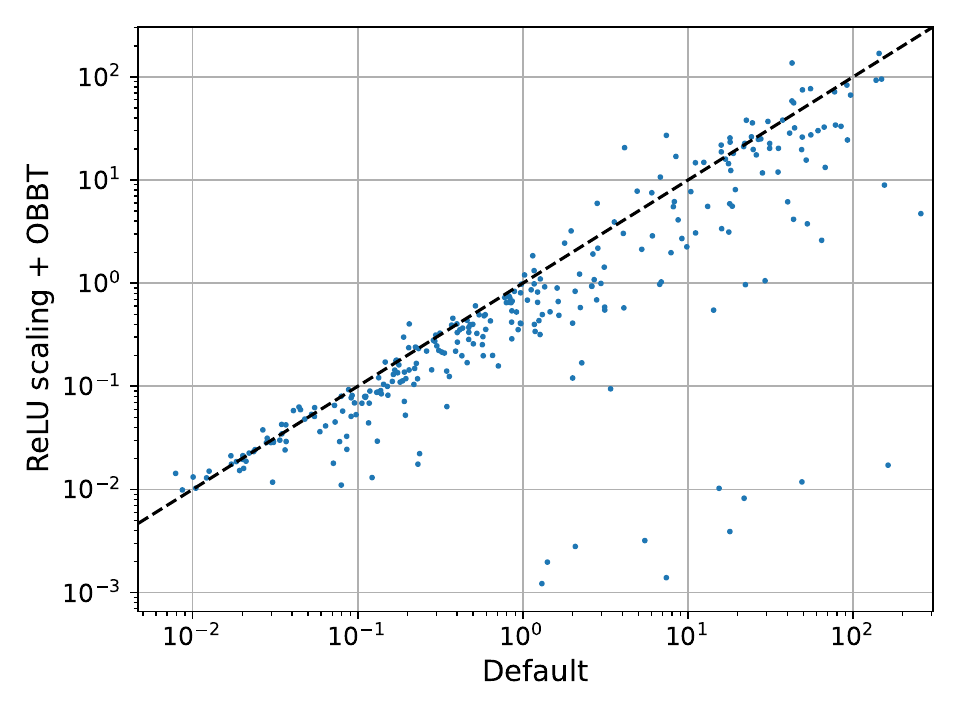}
        \subcaption{Runtime comparison between IA bounds for the baseline network (\enquote{Default}) and OBBT for the scaled ANN (\enquote{ReLU scaling + OBBT}).}
        \label{subfig:IA_vs_scaler_and_OBBT}
    \end{subfigure}

    \caption{Parity plots comparing computational times for optimally solved instances of \eqref{prob:minANN} in different versions: \subref{subfig:IA_vs_scaler}: IA bounds for baseline vs. scaled ReLU network with a geometric mean ratio of {0.936}; \subref{subfig:IA_vs_scaler_and_OBBT}: IA bounds for baseline network vs. OBBT bounds for scaled network with a geometric mean ratio of {0.467}.}
    \label{fig:improvement_scaling_obbt}
\end{figure}

Using the big-M formulation with standard bounds obtained via IA as a baseline, we compare the following options:
\begin{enumerate}
    \item ReLU scaling only: We solve Problem~\eqref{prob:scaling} to obtain equivalent weights and biases with lower $\ell^1$ norm;
    \item ReLU scaling and subsequent LP-based bound tightening: a combination of the two methods.
\end{enumerate}
As shown in \vref{tab:results}, ReLU scaling on its own, as well as combined with OBBT, is able to reduce the big-M coefficients more than applying OBBT on an unscaled network. This is clearly illustrated by the geometric means over the ratios of averaged distances of upper and lower bounds $L$ and $U$ of 0.388 and 0.16 for ReLU scaling and ReLU scaling combined with OBBT compared to unscaled networks, respectively. 
Again, we compute the improvement in computational times as a geometric mean over the ratios of computational times with improved bounds and those with interval arithmetic bounds. We observe that scaling the neural network weights by solving \eqref{prob:scaling} yields only a marginal improvement with a geometric mean of {0.936}. However, combining this scaling with subsequent bound tightening yields a more substantial computational speedup as indicated by a geometric mean of {0.467}. This seems to stem from the tighter big-M coefficients, but also from an increased percentage of stable neurons. Compared to the default bounds, there is an average increase by {7.2} percent. In \vref{fig:improvement_scaling_obbt}, the parities of  computational times for the two comparisons are shown. We note that the parity plot in \vref{subfig:IA_vs_scaler_and_OBBT} suggests that the average speedup may be driven by a few outlier instances in which in the combined method performs exceptionally well.

Overall, with the scaled ReLU networks and their default bounds from interval arithmetic, 307 instances can be solved within the time limit. With tightened bounds, there is a slight reduction to 299 instances.

\subsection{Effect of regularization}

\begin{table}
    \addtolength{\tabcolsep}{-0.2em}
    \centering
    \caption{Influence of training options, bound tightening and ReLU scaling on all trained neural networks and their optimization problems \eqref{prob:minANN}. In each row, the effect of the listed method is evaluated by comparing it to similar networks that differ only in this particular method, e.g., for $\ell^1$ regularization we compare neural networks that were trained with the specified level of regularization to those that were trained without regularization. The first and second column show the number of solved instances without and with the applied technique and the number of instances in total in this comparison. The third column lists the reduction of big-M coefficients as measured by the geometric mean of the ratio of averaged distances of pre-activation bounds $U- L$ of the adapted network and that of the baseline network. The fourth column shows the arithmetic mean of the increase in percentage points of stable neurons due to the applied method. The fifth column shows the geometric mean of the ratio between the number of linear regions of the adapted network and that of the baseline network. The last column shows the geometric mean of the ratio between the computational time with the adapted network and that observed with the baseline network, but is limited to instances in which the optimization problems for both networks are solved within the time limit. We observe a computational speedup with regularization, bound tightening and ReLU-scaling, while dropout leads to a deterioration in performance.}
    \label{tab:results}
    \begin{tabular}{lc|cccccc}
    \multicolumn{1}{c}{}    &  & \begin{tabular}[c]{@{}c@{}}Solved instances\\ (adapted vs. baseline)\end{tabular} & \begin{tabular}[c]{@{}c@{}}Instances\\  total\end{tabular} & 
    \begin{tabular}[c]{@{}c@{}}Geom. mean \\ $\overline{U-L}$ \end{tabular} &
    \begin{tabular}[c]{@{}c@{}}Improvement \\ stable neurons\end{tabular} &    \begin{tabular}[c]{@{}c@{}}Geom. mean \\ lin. regions\end{tabular} & \begin{tabular}[c]{@{}c@{}}Geom. mean \\ time\end{tabular} \\ \hline
    \multirow{5}{*}{$\lambda$} & 1e-3   & 349 vs. 151 & 360 & 0.009  & 0.379    & 0.283  & 0.028 \\
                             & 1e-4   & 358 vs. 151 & 360  & 0.024  & 0.216  & 0.463 & 0.059  \\
                             & 1e-5   & 352 vs. 151 & 360  & 0.052  & 0.122  & 0.817 & 0.109  \\
                             & 1e-6   & 326 vs. 151  & 360 & 0.133  & 0.146  & 1.089 & 0.208  \\
                             & 1e-7   & 287 vs. 151  & 360 & 0.261  & 0.213  & 0.996 & 0.280  \\ \hline
    \multirow{2}{*}{\begin{tabular}[c]{@{}l@{}}Clipped\\ ReLU\end{tabular}} 
                            & M=2    &609 vs. 608 & 720   & 0.415   & 0.029 & 1.094 & 0.931  \\
                            & M=5    & 606 vs. 608 & 720  & 0.560   & 0.011 & 1.055 & 0.974  \\ \hline
    \multirow{2}{*}{Dropout}    & 10\%   & 146 vs. 217 & 240 & 12.761  & -0.204 & 4.098 & 5.825  \\
                                & 20\%   & 152 vs. 217& 240  & 15.403  & -0.210 & 3.513 & 4.845  \\ \hline \hline
    \multicolumn{2}{l|}{OBBT}   & 912 vs. 911  & 1080  & 0.541 & 0.055  & 1.0   & 0.570  \\ \hline
    \multirow{2}{*}{\begin{tabular}[c]{@{}l@{}}ReLU\\ scaling\end{tabular}} &        & 307 vs. 303  & 360   & 0.388 & 0.0    & 1.0   & 0.936  \\
                & OBBT & 299 vs. 303  & 360  & 0.160  & 0.072  & 1.0   & 0.467                                                     
    \end{tabular}
\end{table}

\begin{figure}[H]
    \centering
    \begin{subfigure}{.32\linewidth}
        \centering
        Peaks
        \includegraphics[width=\linewidth]{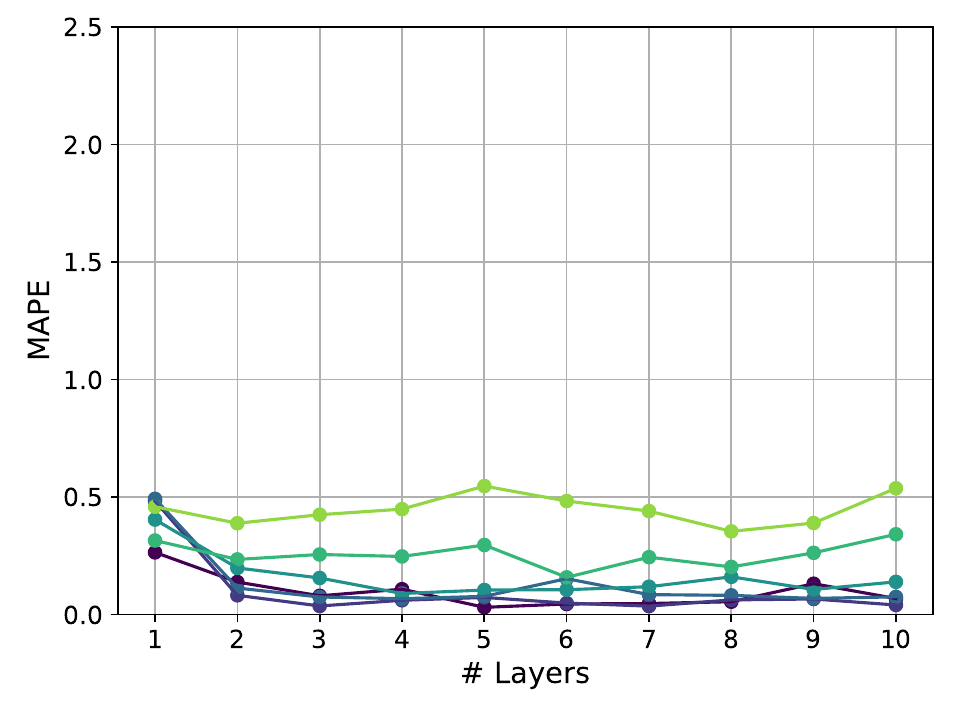}
        \subcaption{MAPE on test set of ReLU ANNs approximating \eqref{fun:peaks}.}
    \end{subfigure}
    \begin{subfigure}{.32\linewidth}
        \centering
        Ackley
        \includegraphics[width=\linewidth]{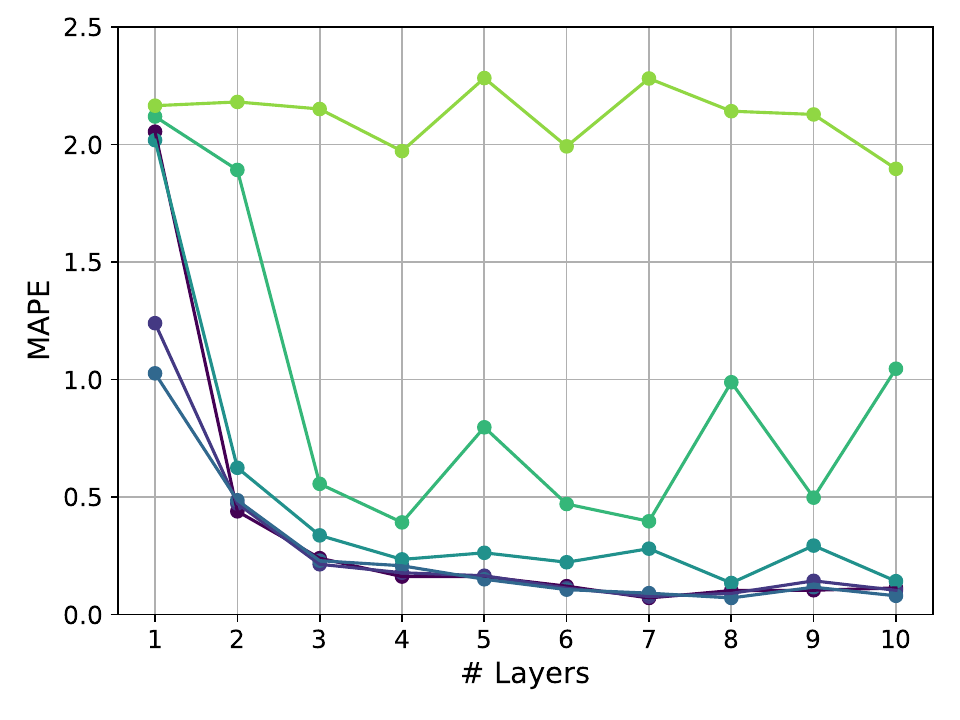}
        \subcaption{MAPE on test set of ReLU ANNs approximating \eqref{fun:ack}.}
    \end{subfigure}
    \begin{subfigure}{.32\linewidth}
        \centering
        Himmelblau
        \includegraphics[width=\linewidth]{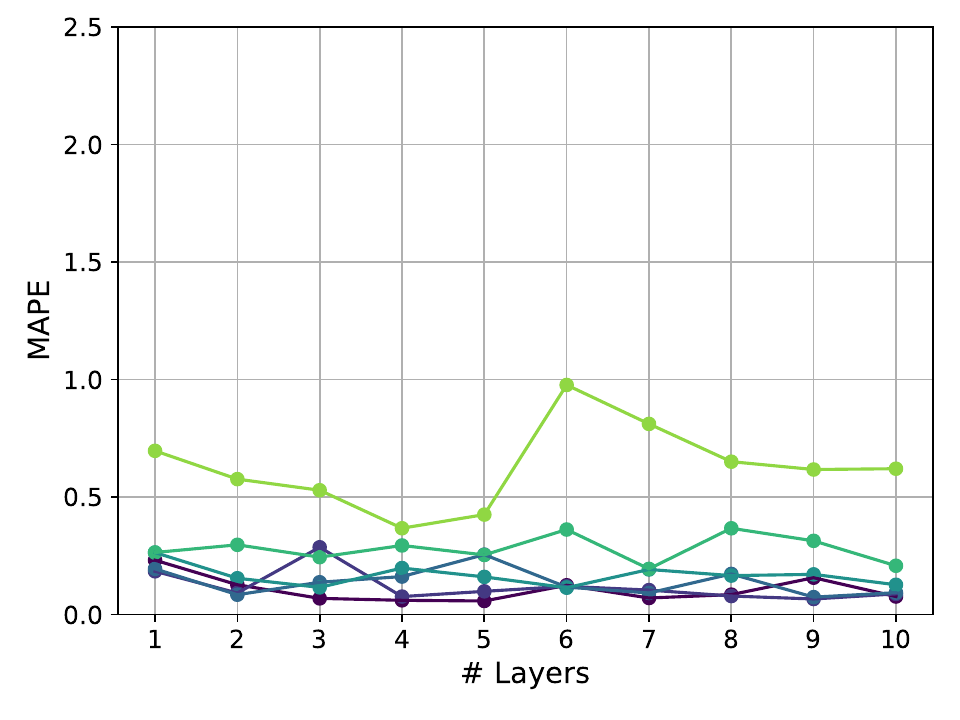}
        \subcaption{MAPE on test set of ReLU ANNs approximating \eqref{fun:him}.}
    \end{subfigure}
    

    \begin{subfigure}{.32\linewidth}
        \centering
        \includegraphics[width=\linewidth]{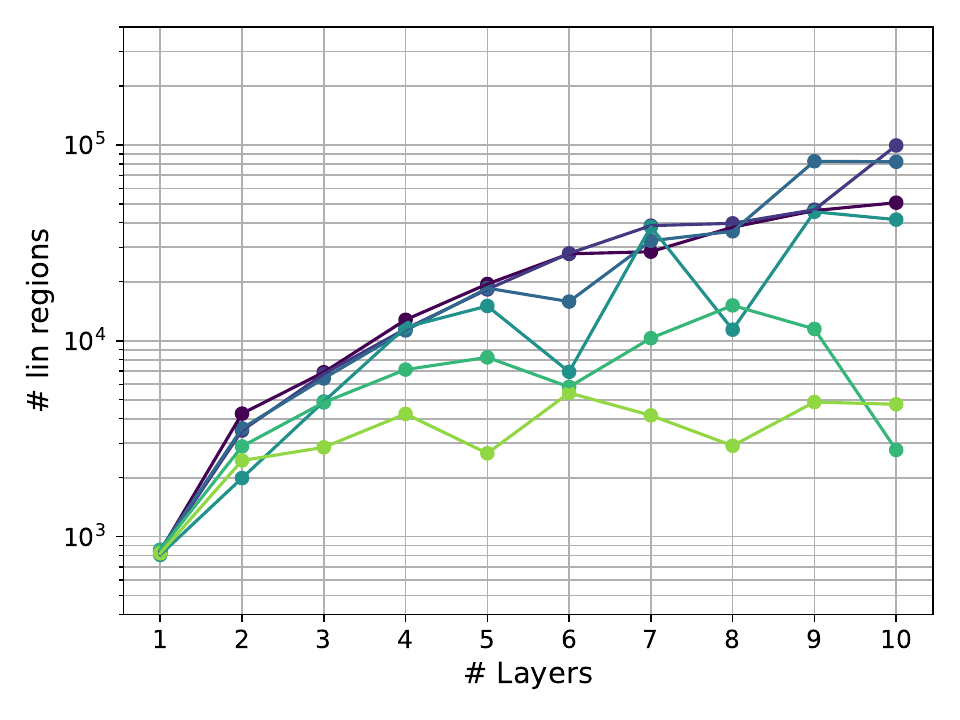}
        \subcaption{Number of linear regions of ReLU ANNs approximating \eqref{fun:peaks}.}
    \end{subfigure}
    \begin{subfigure}{.32\linewidth}
        \centering
        \includegraphics[width=\linewidth]{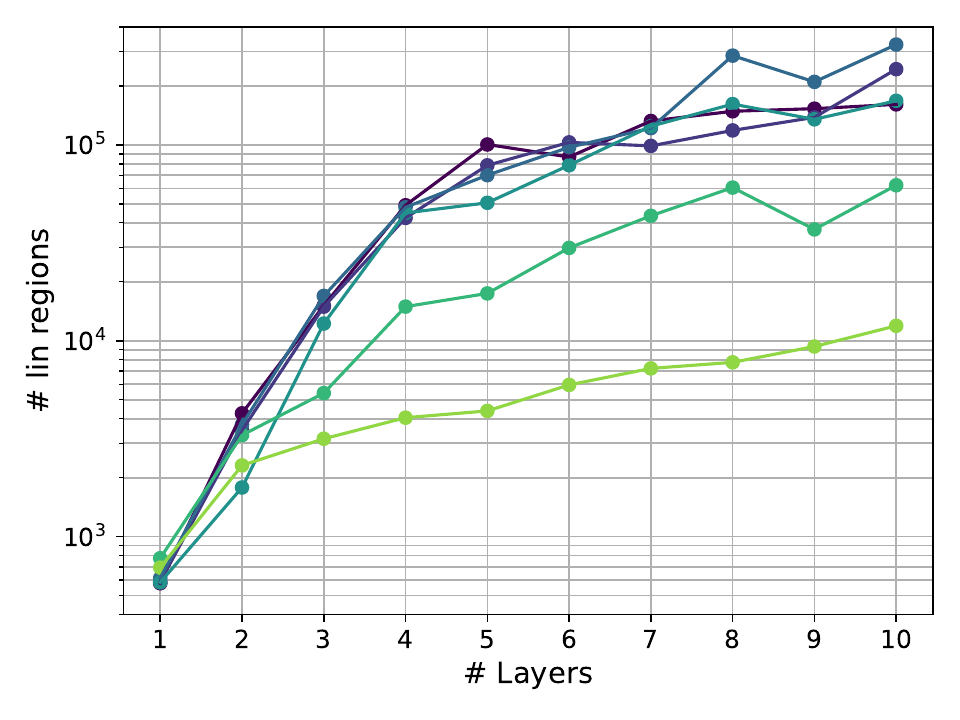}
        \subcaption{Number of linear regions of ReLU ANNs approximating \eqref{fun:ack}.}
    \end{subfigure}
    \begin{subfigure}{.32\linewidth}
        \centering
        \includegraphics[width=\linewidth]{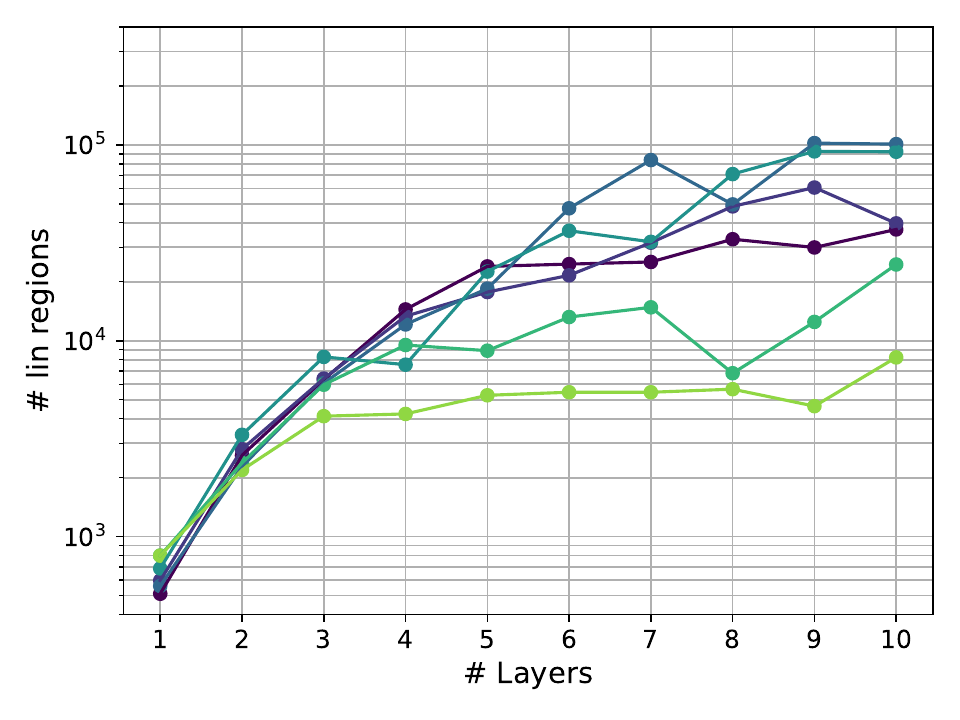}
        \subcaption{Number of linear regions of ReLU ANNs approximating \eqref{fun:him}.}
    \end{subfigure}

    \begin{subfigure}{.32\linewidth}
        \includegraphics[width=\linewidth]{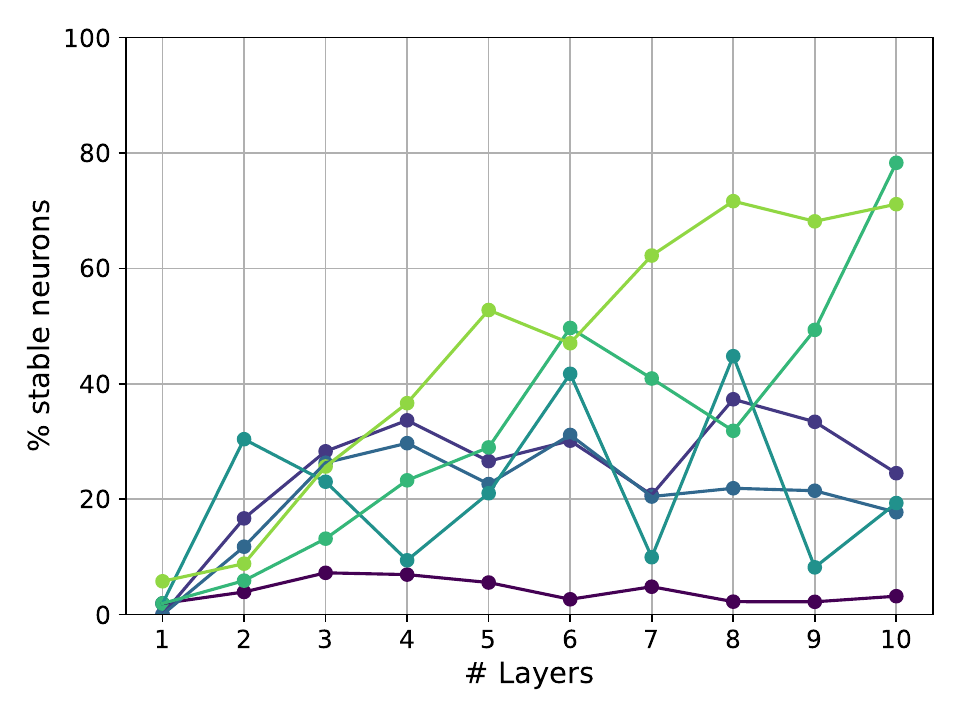}
        \subcaption{Percentage of stable neurons for big-M coefficients based on interval arithmetic bounds.}
    \end{subfigure}
    \begin{subfigure}{.32\linewidth}
        \includegraphics[width=\linewidth]{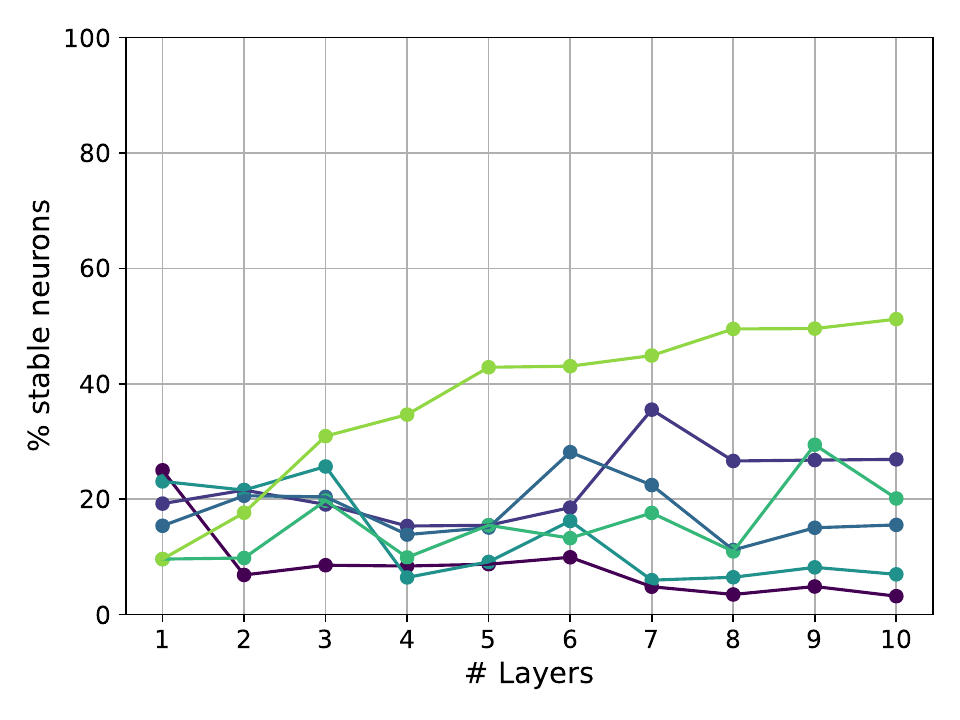}
        \subcaption{Percentage of stable neurons for big-M coefficients based on interval arithmetic bounds.}
    \end{subfigure}
    \begin{subfigure}{.32\linewidth}
        \includegraphics[width=\linewidth]{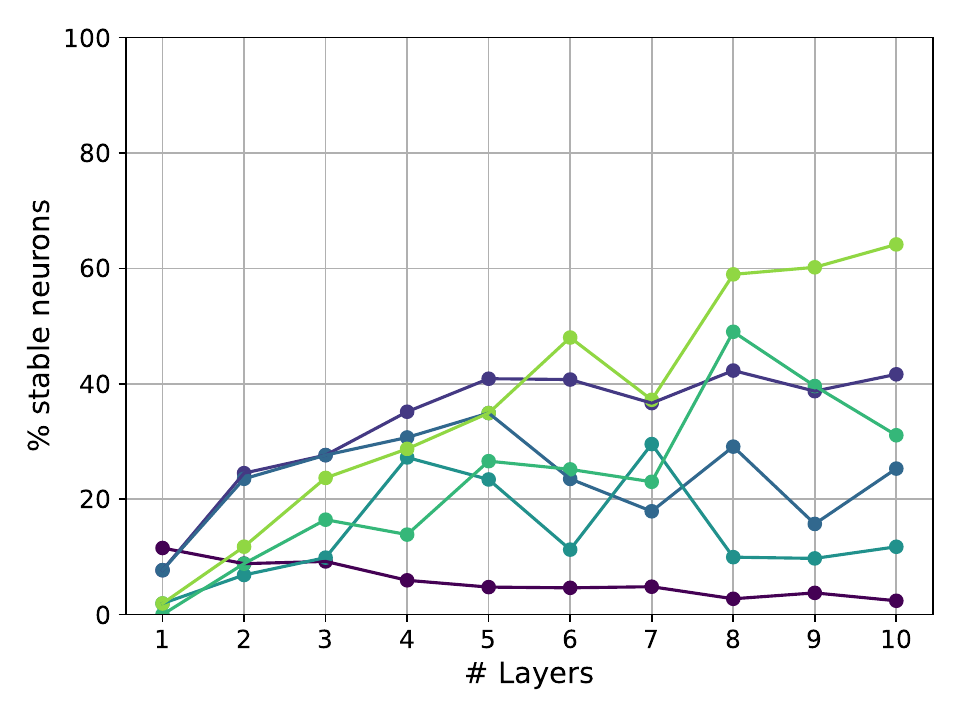}
        \subcaption{Percentage of stable neurons for big-M coefficients based on interval arithmetic bounds.}
    \end{subfigure}

    \begin{subfigure}{.32\linewidth}
        \includegraphics[width=\linewidth]{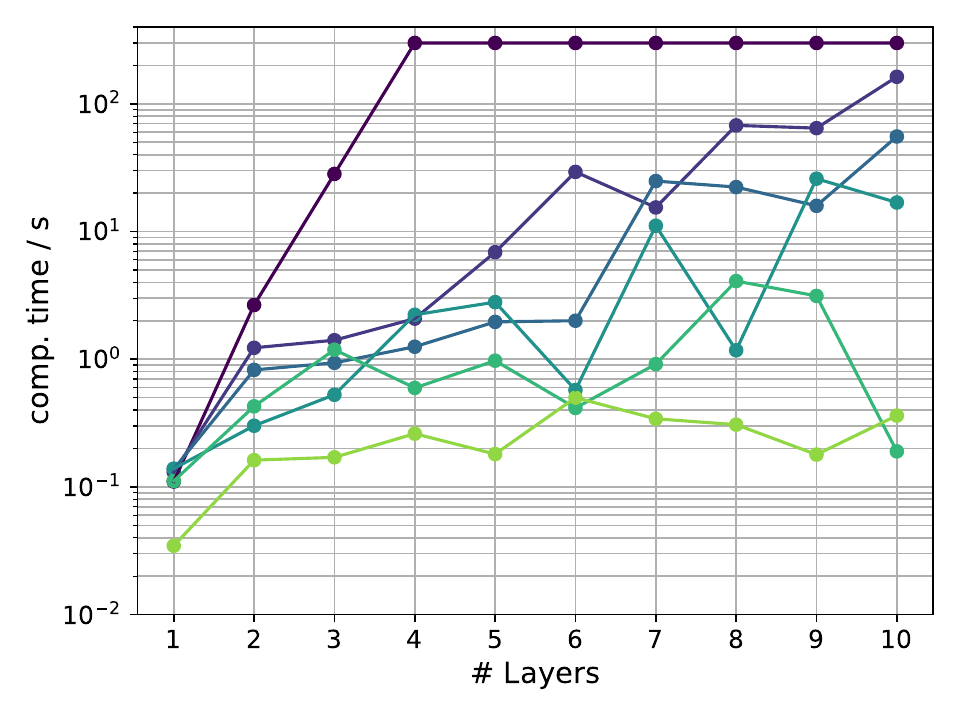}
        \subcaption{Solution time of problem \eqref{prob:minANN} for Peaks function.}
        \typeout{PLOT LINE WIDTH: \the\linewidth}%
    \end{subfigure}
    \begin{subfigure}{.32\linewidth}
        \includegraphics[width=\linewidth]{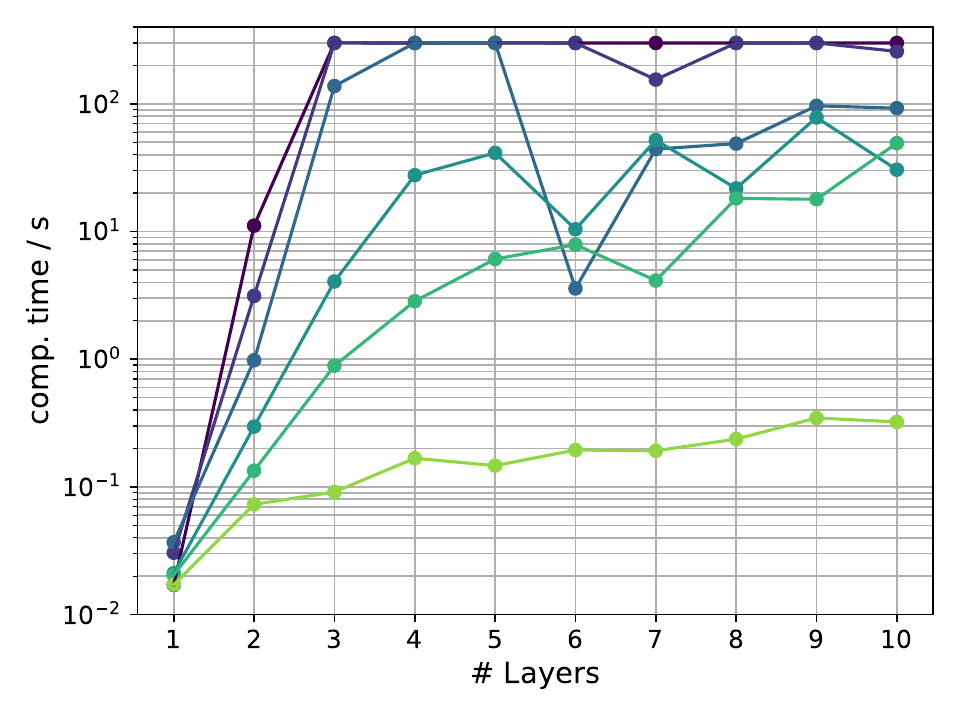}
        \subcaption{Solution time of problem \eqref{prob:minANN} for Ackley's function.}
    \end{subfigure}
    \begin{subfigure}{.32\linewidth}
        \includegraphics[width=\linewidth]{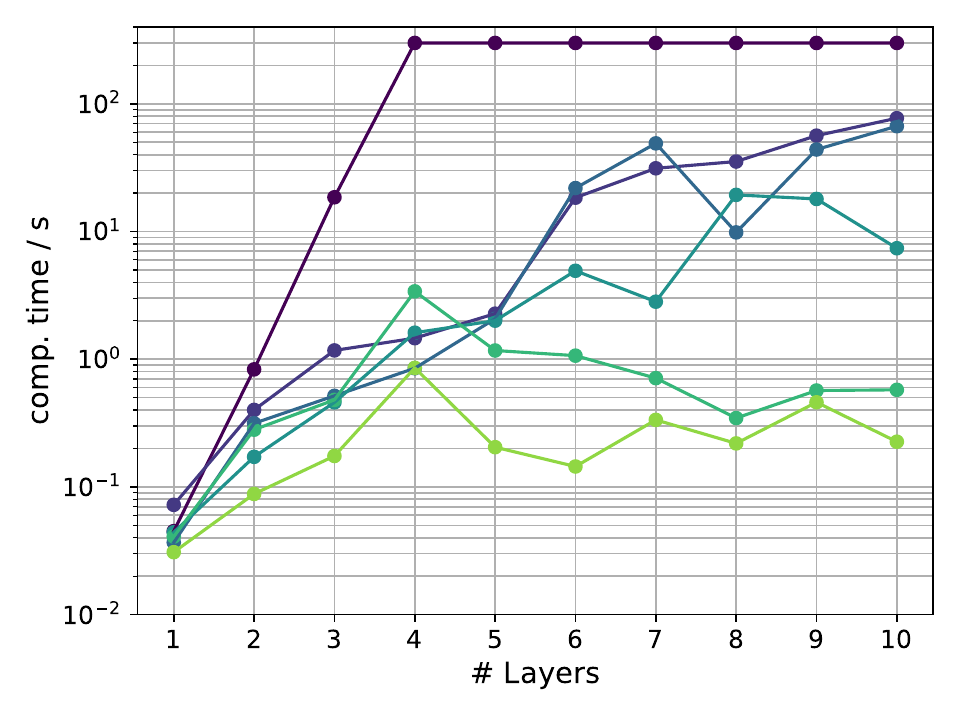}
        \subcaption{Solution time of problem \eqref{prob:minANN} for Himmelblau function.}
    \end{subfigure}
    \begin{subfigure}{\linewidth}
            \centering
            \definecolor{clr1}{HTML}{440154}
\definecolor{clr2}{HTML}{443983}
\definecolor{clr3}{HTML}{31688e}
\definecolor{clr4}{HTML}{21918c}
\definecolor{clr5}{HTML}{35b779}
\definecolor{clr6}{HTML}{90d743}

\begin{tikzpicture} 
    \begin{axis}[%
    hide axis,
    clip bounding box=upper bound,
    xmin=10,
    xmax=50,
    ymin=0,
    ymax=0.4,
    line width= 2,
    scale only axis,
    no markers,
    enlargelimits=false,
    legend columns = 7,
    legend style={draw=none,legend cell align=left}
    ]
    \addlegendimage{empty legend}
    \addlegendentry{$\ell^1$ regularization:}
    \addlegendimage{clr1}
    \addlegendentry{0};
    \addlegendimage{clr2}
    \addlegendentry{1e-7};
    \addlegendimage{clr3}
    \addlegendentry{1e-6};
    \addlegendimage{clr4}
    \addlegendentry{1e-5};
    \addlegendimage{clr5}
    \addlegendentry{1e-4};
    \addlegendimage{clr6}
    \addlegendentry{1e-3};

\end{axis}
\end{tikzpicture}
    \end{subfigure}
    \caption{Mean absolute percentage error on test set, number of linear regions, percentage of stable neurons and computational times in problem \eqref{prob:minANN} of trained ANNs with varying number of hidden layers with 50 neurons, trained with different levels of $\ell^1$ regularization.}
    \label{fig:model_statistics}%
    \typeout{LINE WIDTH: \the\linewidth}%
\end{figure}

\begin{figure}
    \centering
    \begin{subfigure}{.32\linewidth}
        \centering
        \includegraphics[width=\linewidth,trim={0 0 0 8mm},clip]{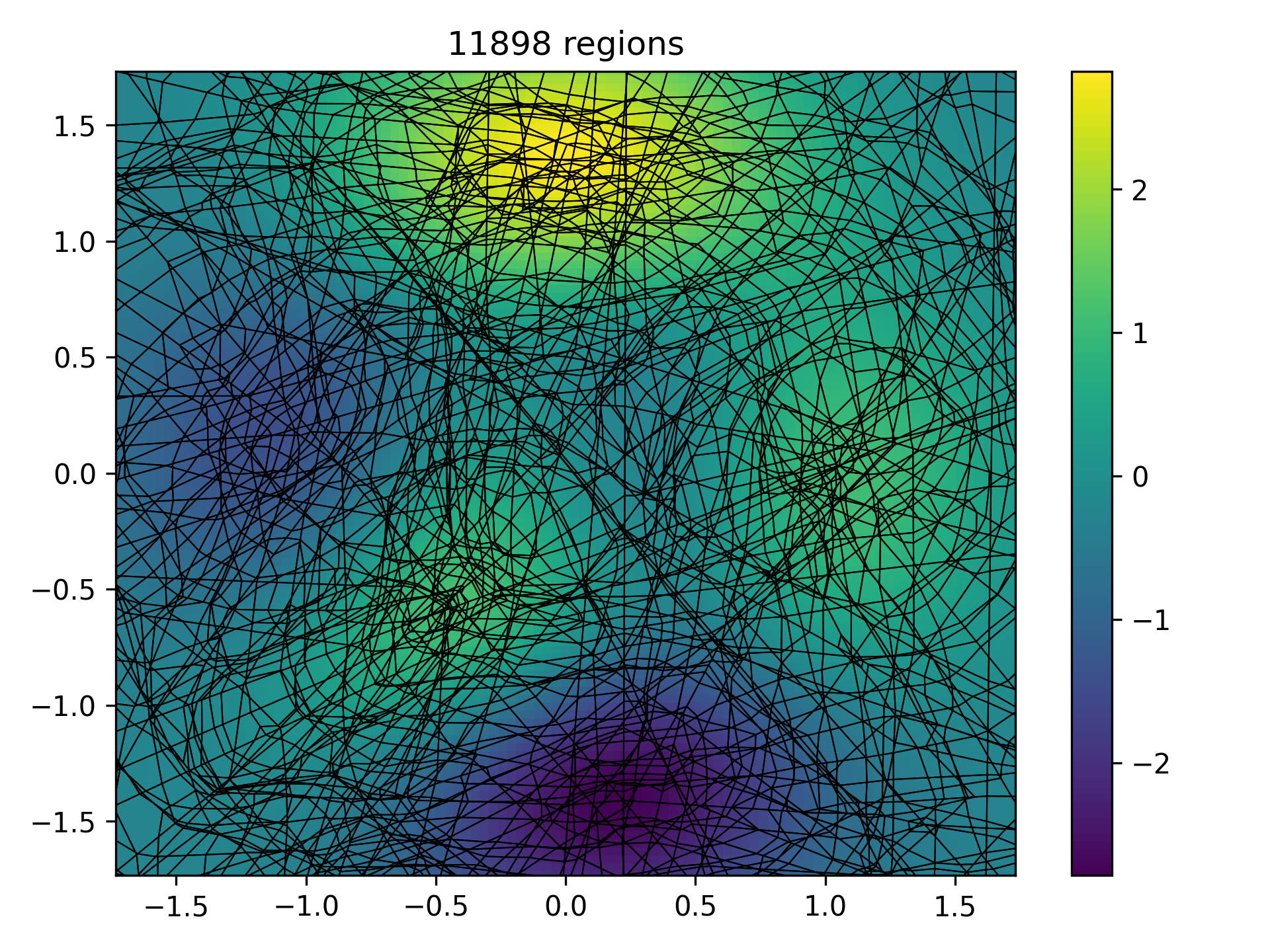}
        \subcaption{Standard, 11,898 regions.}
        \label{subfig:mesh_no_regu}
    \end{subfigure}
    \begin{subfigure}{.32\linewidth}
        \centering
        \includegraphics[width=\linewidth,trim={0 0 0 8mm},clip]{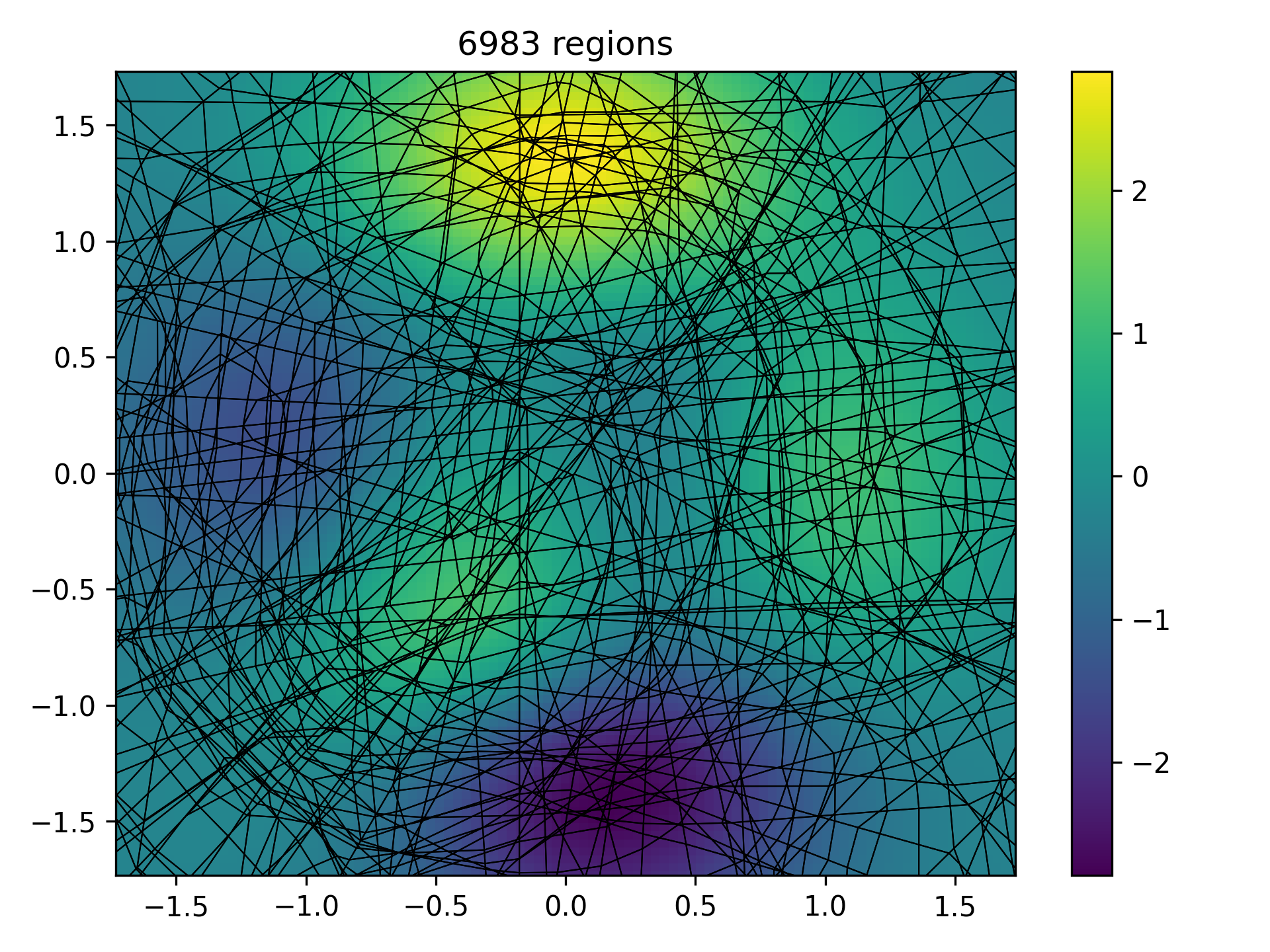}
        \subcaption{Regularization, 6,893 regions.}
        \label{subfig:mesh_regu}
    \end{subfigure}
    \begin{subfigure}{.32\linewidth}
        \centering
        \includegraphics[width=\linewidth,trim={0 0 0 8mm},clip]{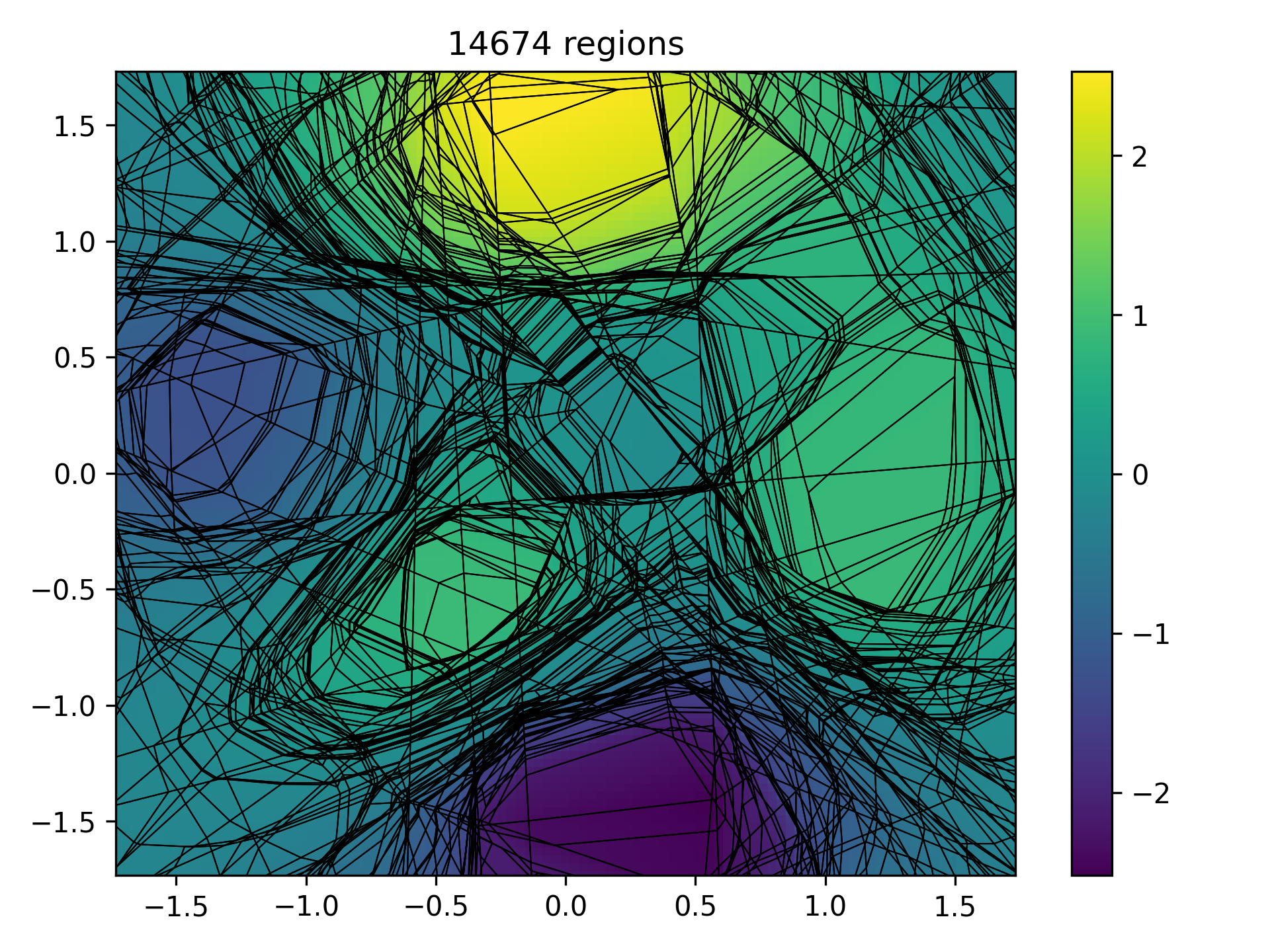}
        \subcaption{Dropout, 14,674 regions.}
        \label{subfig:mesh_dropout}
    \end{subfigure}
    \caption{Linear regions for ReLU networks approximating the Peaks function \eqref{fun:peaks} with five hidden layers of 25 neurons each. Color-coded in the backgrounds are the outputs of the neural networks. Compared are networks with different training options: \cref{subfig:mesh_no_regu} with no regularization or dropout, \cref{subfig:mesh_regu} with $\ell^1$ regularization and $\lambda=10^{-5}$, \cref{subfig:mesh_dropout} with $20$\% dropout. Regularizing the weights of the ANN during training decreases the number of linear regions, applying dropout increases it and also changes their sizes.}
    \label{fig:linear_regions}
\end{figure}

In \vref{fig:model_statistics}, we depict how the mean absolute error on the test set, the number of linear regions, the percentage of neurons of fixed activation, and the solver runtime correlate with the depth of the neural networks for networks with 50 neurons per layer with different regularization parameters. 
In the first row, we see the performance of the neural networks on the test data as measured by the MAPE. We see, that large regularization parameters lead to a degradation of accuracy on the test dataset, especially for Ackley's function. For small regularization parameters there is a high level of agreement between the predictions and the ground truth on the test data. Further, in some instances, training with moderate levels of $\ell^1$ regularization does in fact lead to be better generalization of the neural network. 
The second row of \cref{fig:model_statistics} shows the number of linear regions as an indicator of the complexity, or expressive power of the neural network. With increasing levels of regularization, we obtain neural networks with a lower number of linear regions. This is also illustrated in \vref{fig:linear_regions}. Comparing the number of linear regions among the three different functions, the neural networks which approximate Ackley's function have the most linear regions. This is plausible comparing the surface plots in \vref{fig:different_functions}, because Ackley's function exhibits a large number of local oscillations.
In the third row of \cref{fig:model_statistics}, we plot the percentage of stable neurons. These are neurons whose input bounds are either non-negative or non-positive, which means that they are in a fixed state of activation regardless of input. No binary variables have to be added to model the activation function of such neurons. Confirming the findings of \citet{Xiao2019,Serra2020}, higher values of $\lambda$ lead to a higher percentage of stable neurons. 
The last row shows the computational times in the optimization problem \eqref{prob:minANN}. Comparing the runtimes among the three functions, Ackley's function appears to be the hardest to minimize. Here, we cannot solve unregularized networks with as little as three hidden layers to global optimality within the specified time limit. Based on the observation that ANNs approximating this function have an increased number of linear regions and that several local minima exist in the input domain, this is expected behavior. Increasing the regularization generally lowers the time to compute global minima for all three functions. While the global minima of unregularized networks cannot be determined for any network with more than four layers, applying moderate levels of regularization makes almost all instances tractable. The only exception here is Ackley's function, which remains unsolved for the lowest regularization parameter $\lambda = 10^{-7}$ as well. While \cref{fig:model_statistics} shows the results for all networks with 50 neurons per hidden layer, we obtain similar results for those with 25 neurons (data shown in the appendix).
In combination with the results in \vref{tab:results}, this illustrates that regularization proved the most effective method by improving big-M coefficients, increasing the number of stable neurons and decreasing the number of linear regions, thus enabling the computational speedup.

\subsection{Effect of clipped ReLU}

\begin{figure}[ht!]
    \centering
    \begin{subfigure}{0.49\textwidth}
        \includegraphics[width=\textwidth]{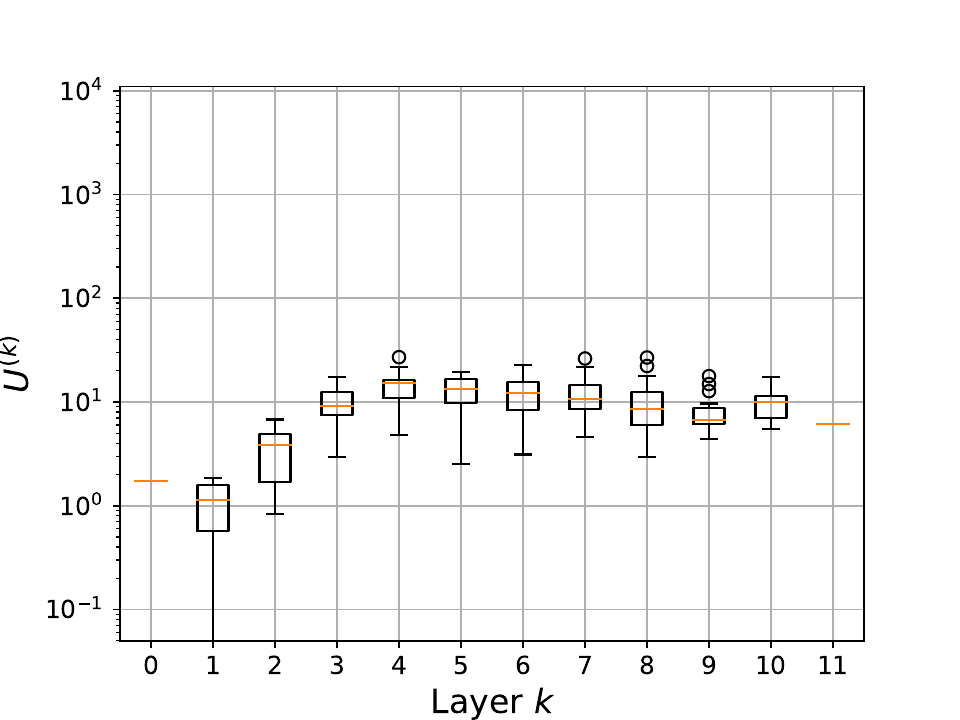}
        \caption{Pre-activation bounds $U^{(k)}$ determined via interval arithmetic for clipped ReLU ANN with the big-M formulation \eqref{eq:bigM_clipped} and $M=5.0$.}
        \label{fig:IA_bounds_clipped}
    \end{subfigure}
    \hfill
    \begin{subfigure}{0.49\textwidth}
        \includegraphics[width=\textwidth]{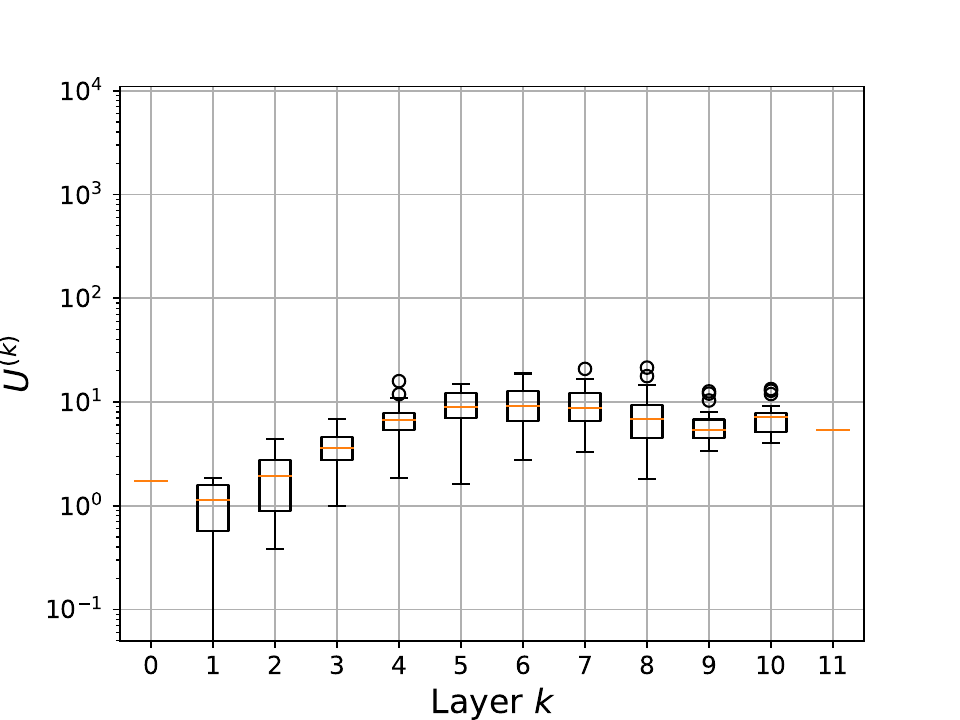}
        \caption{Pre-activation bounds $U^{(k)}$ determined via LP-based bound tightening for clipped ReLU ANN with big-M formulation \eqref{eq:bigM_clipped} and $M=5.0$.}
        \label{fig:LR_bounds_clipped}
    \end{subfigure}
   
    \caption{Comparison of pre-activation bounds  $U^{(k)}$ for neural networks with ten hidden layers and clipped ReLU formulation \eqref{eq:bigM_clipped} with $M=5.0$ as activation function. Compared to the bounds derived via interval arithmetic for the regular ReLU activation shown in \vref{fig:IA_bounds}, the bounds for the clipped ReLU are generally lower. Moreover, due to the threshold $M$, the bounds stay approximately constant over the layers. Solving auxiliary LPs only noticeably tightens bounds in the first few layers.}
    \label{fig:big_Ms_clipped}
\end{figure}

The effect of the clipping is obvious in the big-M coefficients of formulation \eqref{eq:bigM_clipped}, which are illustrated in \vref{fig:big_Ms_clipped} for a threshold of $M=5.0$. Compared to the big-M coefficients derived for the regular ReLU activation function and depicted in \vref{fig:IA_bounds}, the clipped ReLU formulation yields lower bounds, though this may depend on the particular choice of $M$. This is also obvious from the results in \vref{tab:results}, with clipped ReLU leading to greater reductions in big-M coefficients compared to OBBT. We also observe that LP-based bound tightening for neural networks with clipped ReLU activation does not improve the bounds to the same degree as it did for the regular ReLU activation function as depicted in \vref{fig:LR_bounds}.

As the results in \Cref{tab:results} suggest, using the clipped ReLU  \eqref{eq:bigM_clipped} yields only marginal computational speedup compared to the standard ReLU activation. There seems to be a tradeoff between a higher number of binary variables needed for modeling \eqref{eq:bigM_clipped} and the slightly higher number of linear regions on the one hand, and the reduction of big-M coefficients on the other hand. 

\subsection{Effect of Dropout}

For the Peaks function only, we trained additional networks with different levels of dropout applied to the hidden layers, namely 10 and 20 percent. In combination with the other hyperparameters (regularization, depth and width of the network) this yields a total of 240 trained neural networks with dropout whose properties we can compare. As illustrated in \Cref{tab:results}, we find that dropout leads to neural networks with three to four times more linear regions on average, confirming the findings of \citet{Zhang2020a}. This is also evident in \Cref{fig:linear_regions}, where the linear regions of three exemplary ANNs are compared for networks with five hidden layers. Another effect of dropout is the percentage of stable neurons, which is reduced by approximately 20~\% on average compared to networks trained without dropout and a drastic increase in the magnitude of the big-M coefficients. In combination, this leads to a reduction in instances that could be solved to global optimality and a simultaneous four to six-fold increase in computational time for those instances that could be solved.

\section{Conclusions and outlook}\label{sec:discussion}

In this paper, we compared different variations of training and scaling methods for ReLU networks with respect to their effect on the performance of global optimization solvers on problems with these networks  embedded. We divided these methods into those that are applied during training and those that can be used on trained networks. For the latter category, we proposed a scaling method specific to the ReLU activation function, which equivalently transforms a ReLU ANN such that the $\ell^1$ norm of the networks weights and biases is minimized. This has the desired effect of reducing the constant coefficients in big-M formulations of the network's activation functions. In numerical experiments, we demonstrated that this method can be used to reduce the computational effort of solving subsequent optimization problems, when it is used in combination with bound tightening. Although in our study we only investigated the direct minimization of feed-forward neural networks with their big-M formulation of ReLU networks, we believe that the findings are also applicable in other contexts. These might include optimization problems with ReLU networks using different MILP encodings, e.g., the partition-based formulation from \citet{Tsay2021}, or other optimization settings, e.g., more difficult optimization problems from real-world applications. In fact, by employing regularization during training we were able to solve a complex superstructure optimization problem in chemical engineering that had been computationally intractable before \citep{Klimek2024}.

Moreover, to the best of our knowledge, this is the first computational study that links various training methods to both the number of linear regions and the percentage of fixed neurons as well as the computational effort in subsequent optimization problems. Doing so, we were able to provide empirical evidence for several observations from the literature, e.g., an increased number of linear regions for networks trained with dropout, and computational speedup due to higher rates of fixed neurons for networks trained with $\ell^1$ regularization. 

Further research may include a more thorough analysis into how the used training methods and hyperparameter options used when training a neural network impact its number of linear regions and the number of fixed neurons. Also, different objectives in \eqref{prob:scaling} may be conceivable to promote other properties in the transformed networks. It may also be promising to investigate transformations that allow minor perturbations of the functional relationship.

\begin{acknowledgements}
This project has received funding from the European Regional Development Fund (grants timingMatters and IntelAlgen) under the European Union’s Horizon Europe Research and Innovation Program, 
from the research initiative ``SmartProSys: Intelligent Process Systems for the Sustainable Production of Chemicals'' funded by the Ministry for Science, Energy, Climate Protection and the Environment of the State of Saxony-Anhalt, and
from the German Research Foundation DFG within GRK 2297 ’Mathematical Complexity Reduction’ and priority program 2331 ’Machine Learning in Chemical Engineering’ under grant SA 2016/3-1.
\end{acknowledgements}

\bibliographystyle{spbasic}
\bibliography{main}

\begin{appendices}
\section{Results for smaller models}
\begin{figure}[h]
    \centering
    \begin{subfigure}{.32\linewidth}
        \centering
        Peaks
        \includegraphics[width=\linewidth]{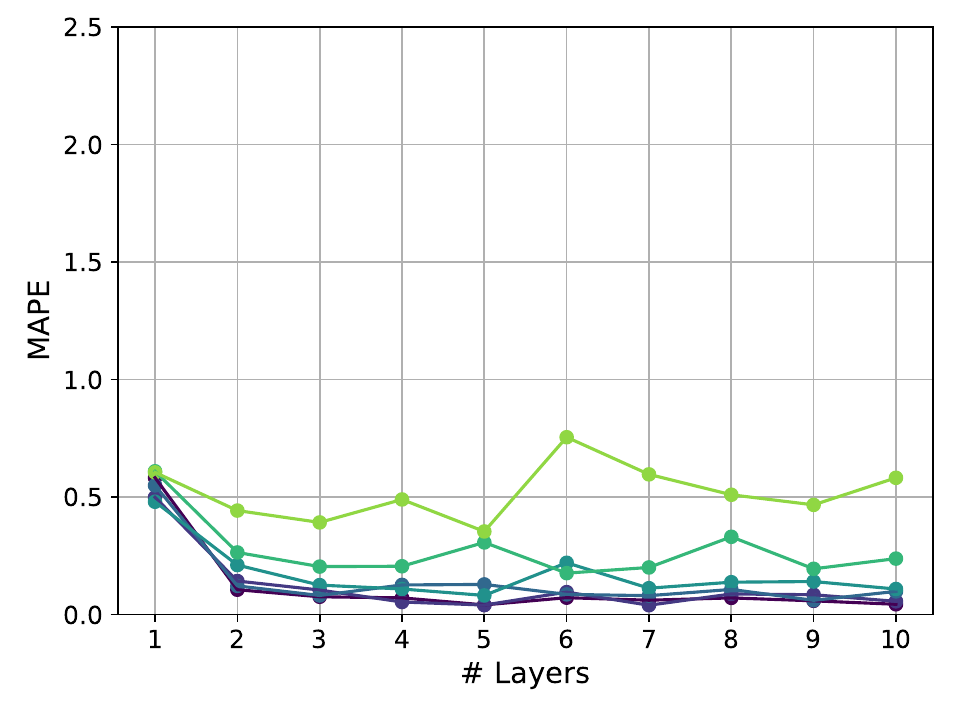}
        \subcaption{MAPE on test set of ReLU ANNs approximating \eqref{fun:peaks}.}
    \end{subfigure}
    \begin{subfigure}{.32\linewidth}
        \centering
        Ackley
        \includegraphics[width=\linewidth]{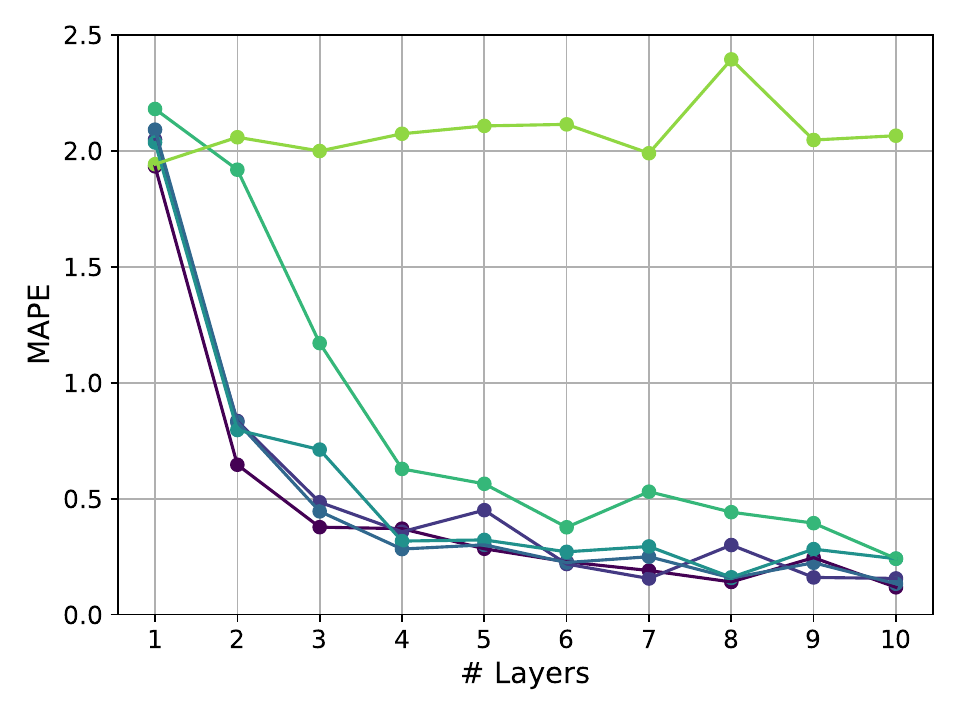}
        \subcaption{MAPE on test set of ReLU ANNs approximating \eqref{fun:ack}.}
    \end{subfigure}
    \begin{subfigure}{.32\linewidth}
        \centering
        Himmelblau
        \includegraphics[width=\linewidth]{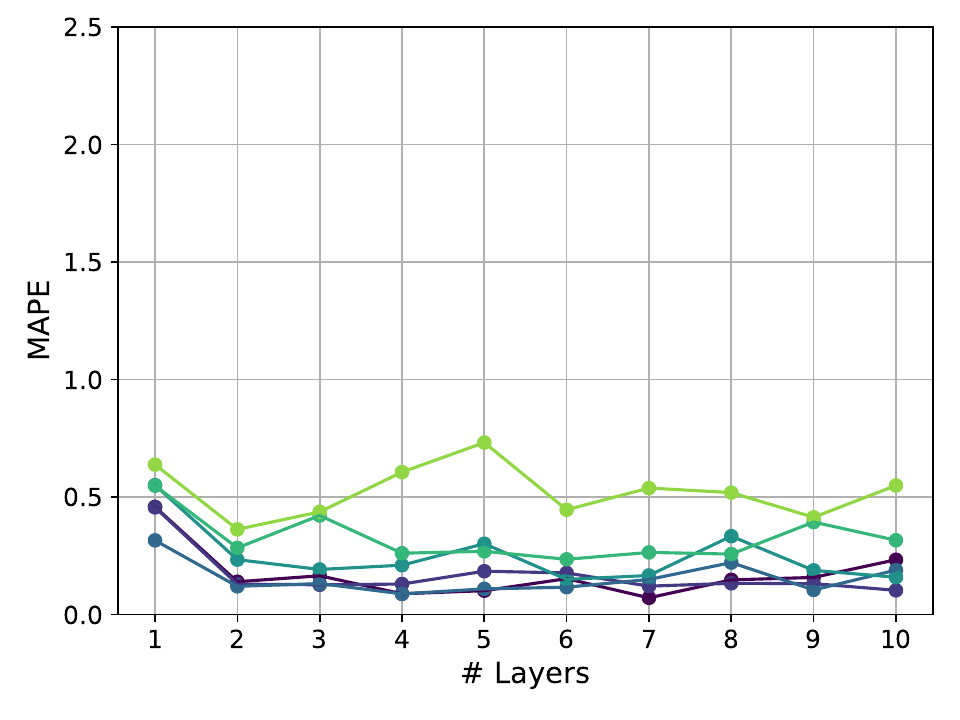}
        \subcaption{MAPE on test set of ReLU ANNs approximating \eqref{fun:him}.}
    \end{subfigure}
    
    \begin{subfigure}{.32\linewidth}
        \centering
        \includegraphics[width=\linewidth]{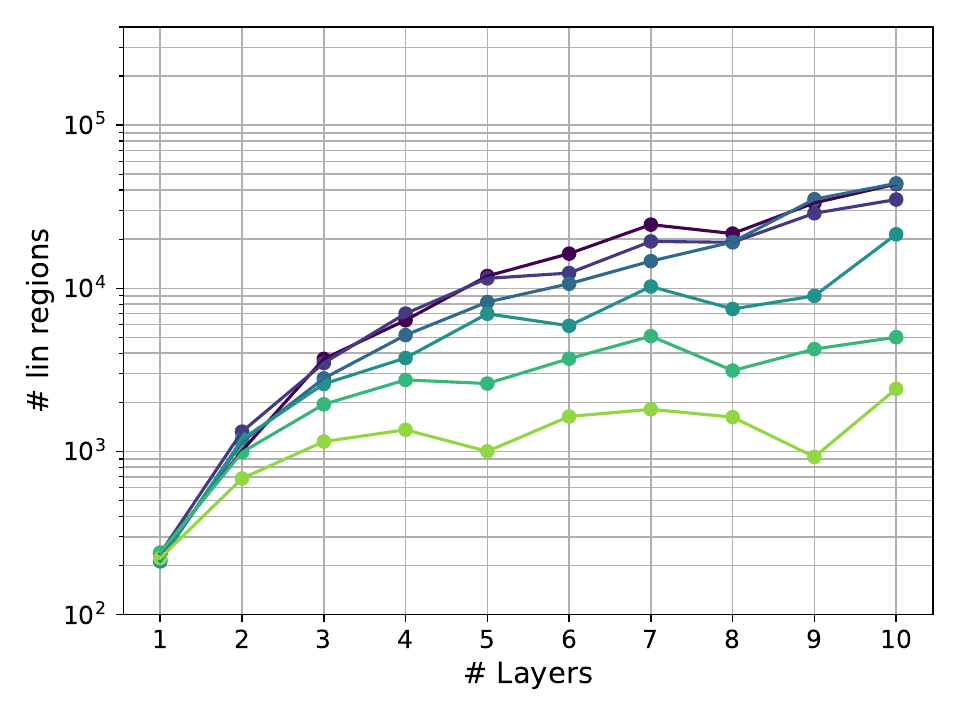}
        \subcaption{Number of linear regions of ReLU ANNs approximating \eqref{fun:peaks}.}
    \end{subfigure}
    \begin{subfigure}{.32\linewidth}
        \centering
        \includegraphics[width=\linewidth]{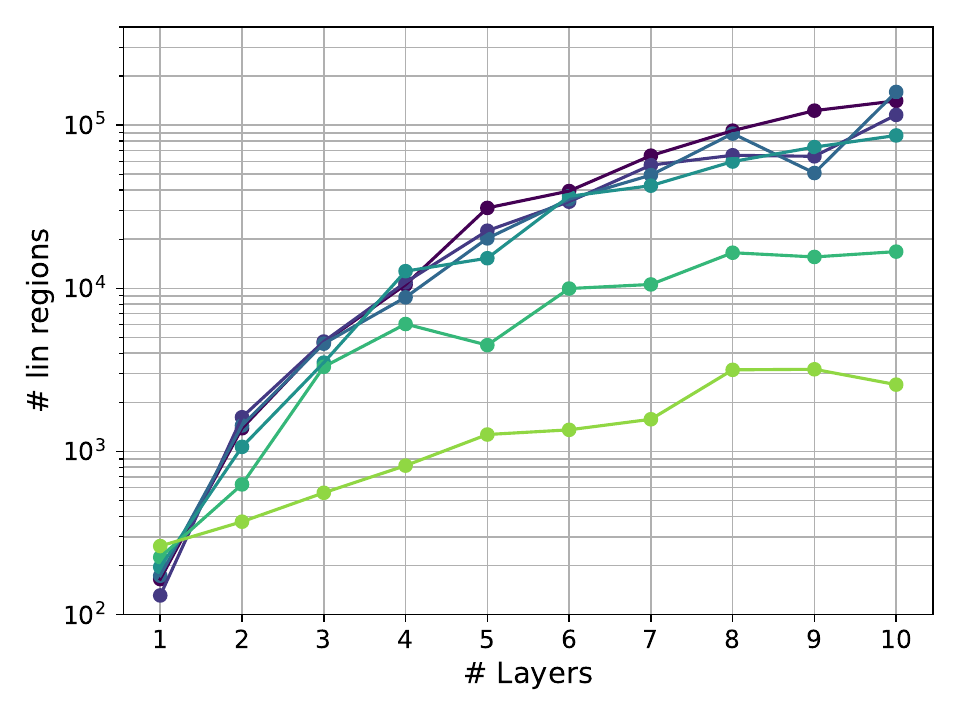}
        \subcaption{Number of linear regions of ReLU ANNs approximating \eqref{fun:ack}.}
    \end{subfigure}
    \begin{subfigure}{.32\linewidth}
        \centering
        \includegraphics[width=\linewidth]{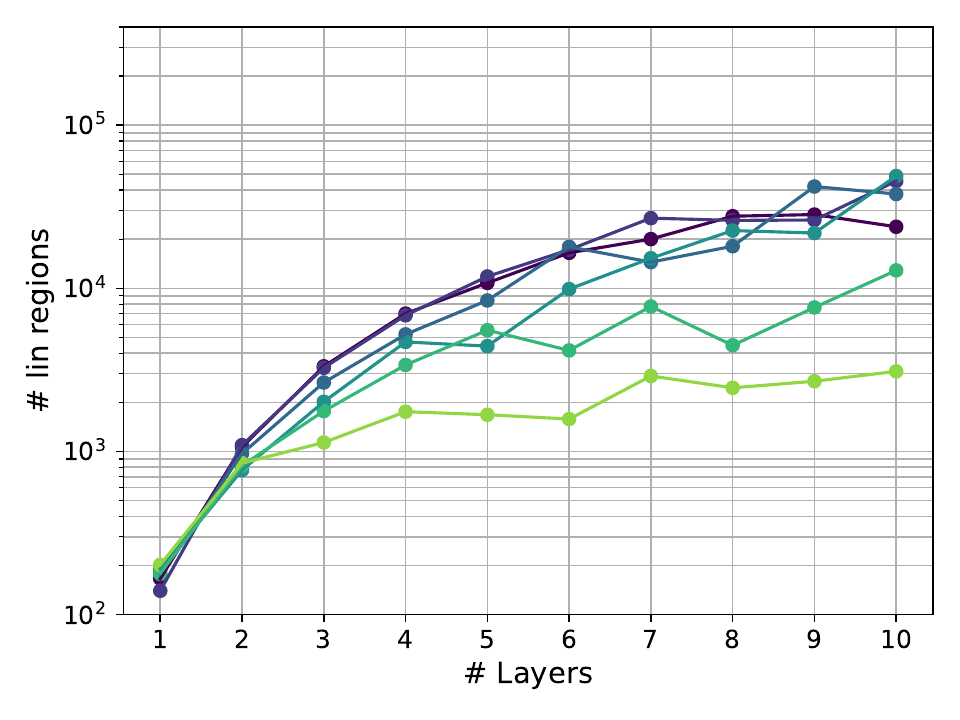}
        \subcaption{Number of linear regions of ReLU ANNs approximating \eqref{fun:him}.}
    \end{subfigure}

    \begin{subfigure}{.32\linewidth}
        \includegraphics[width=\linewidth]{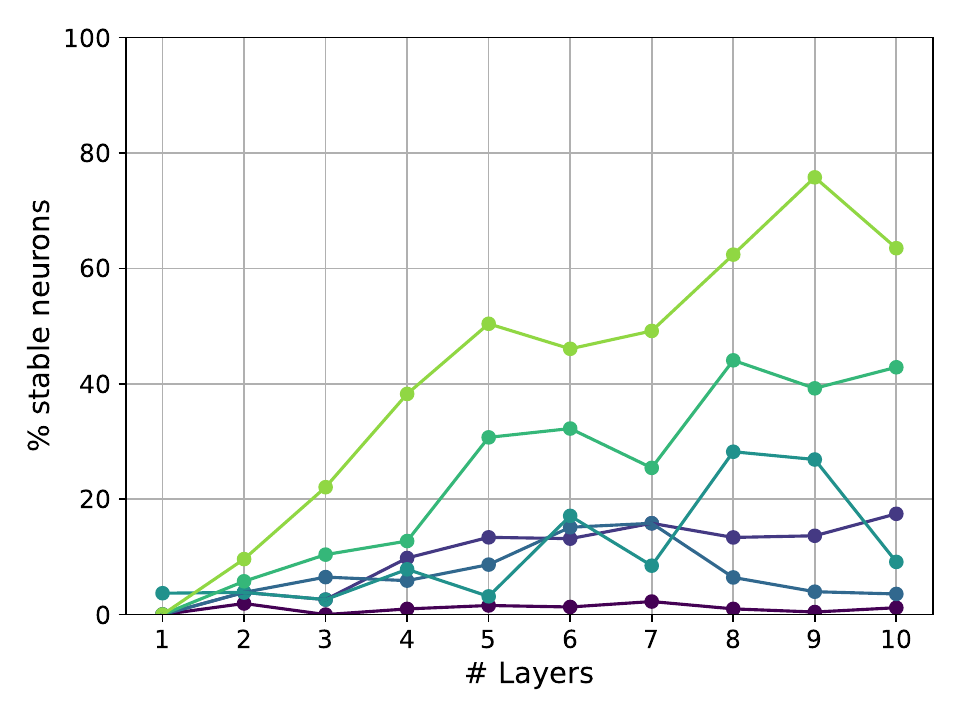}
        \subcaption{Percentage of fixed neurons based on IA bounds.}
    \end{subfigure}
    \begin{subfigure}{.32\linewidth}
        \includegraphics[width=\linewidth]{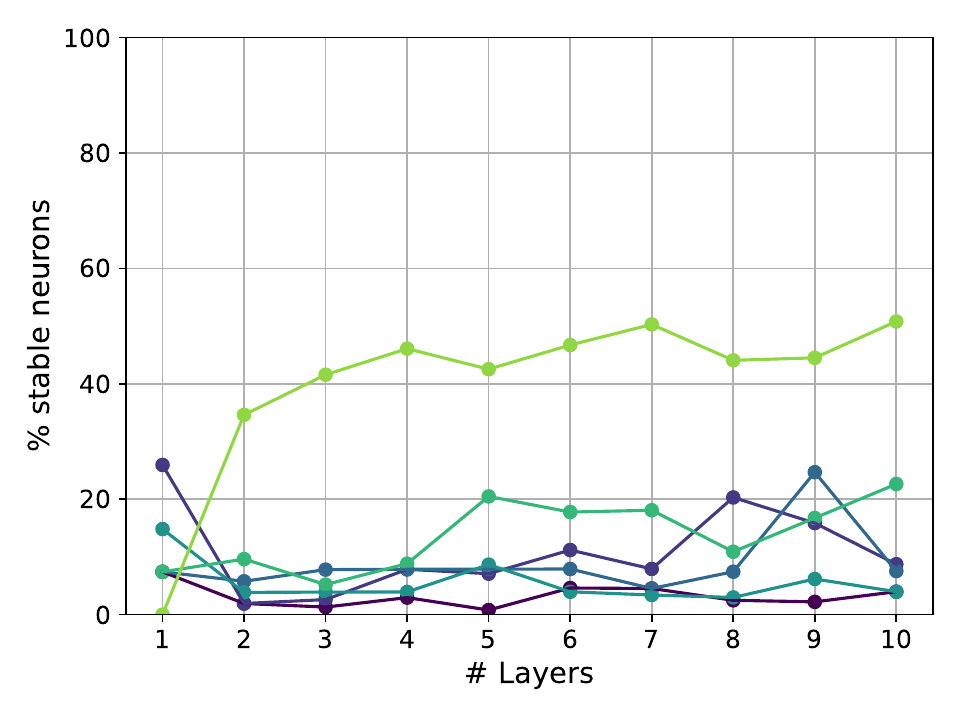}
        \subcaption{Percentage of fixed neurons based on IA bounds.}
    \end{subfigure}
    \begin{subfigure}{.32\linewidth}
        \includegraphics[width=\linewidth]{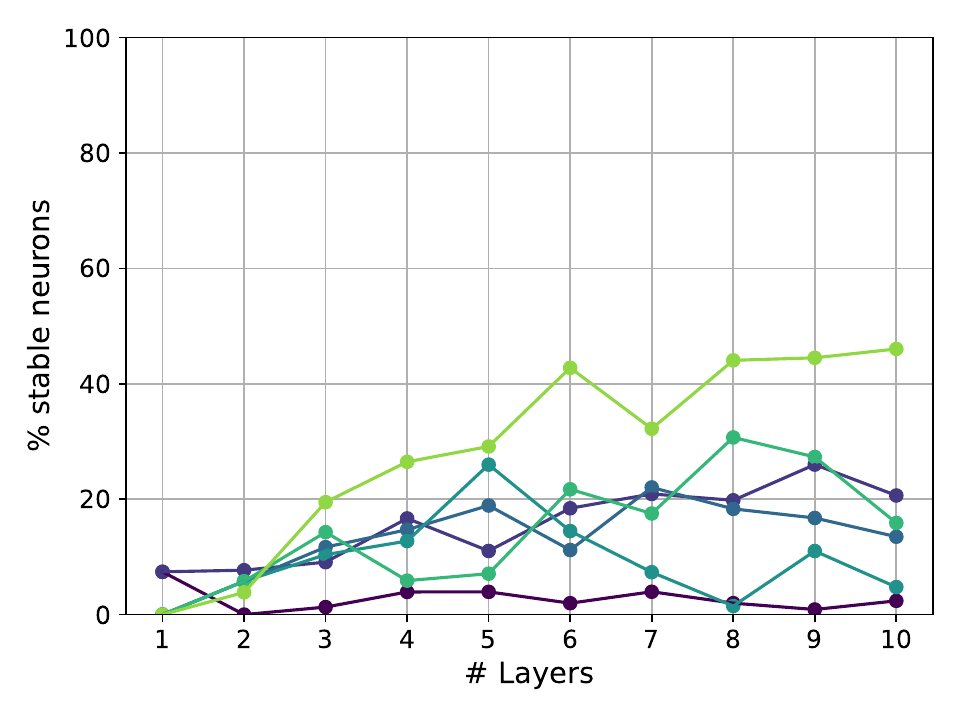}
        \subcaption{Percentage of fixed neurons based on IA bounds.}
    \end{subfigure}

    \begin{subfigure}{.32\linewidth}
        \includegraphics[width=\linewidth]{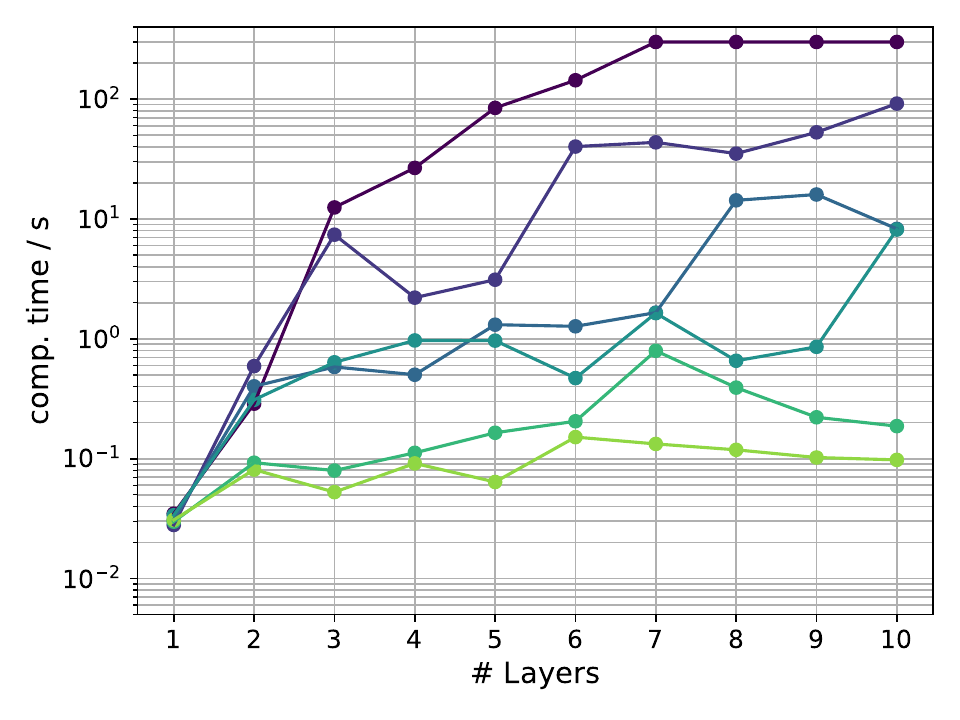}
        \subcaption{Solution time of problem \eqref{prob:minANN} for Peaks function.}
    \end{subfigure}
    \begin{subfigure}{.32\linewidth}
        \includegraphics[width=\linewidth]{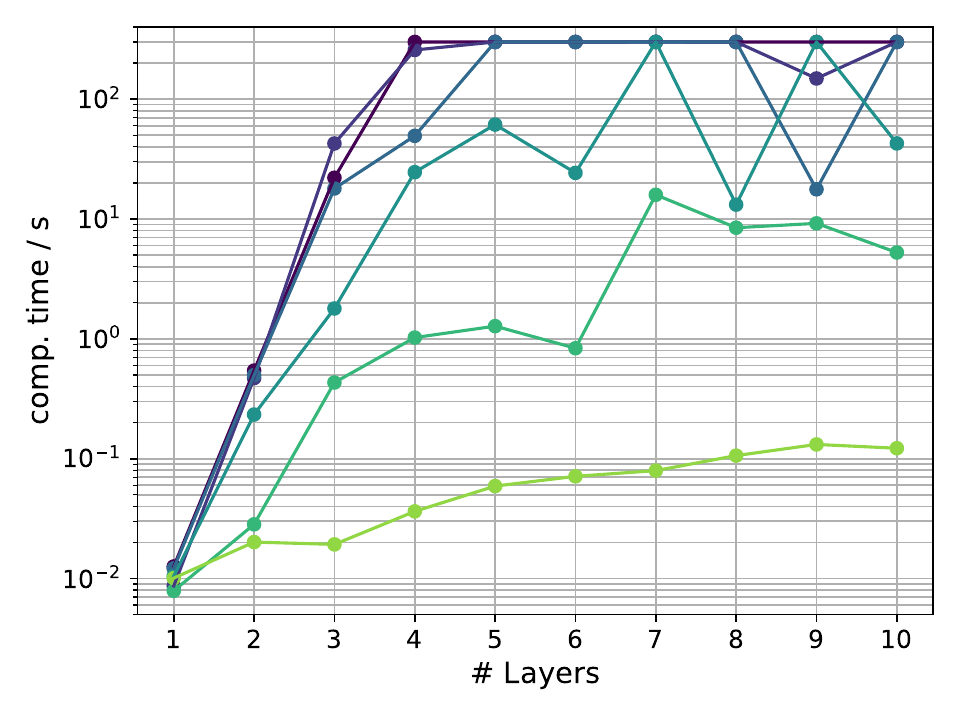}
        \subcaption{Solution time of problem \eqref{prob:minANN} for Ackley's function.}
    \end{subfigure}
    \begin{subfigure}{.32\linewidth}
        \includegraphics[width=\linewidth]{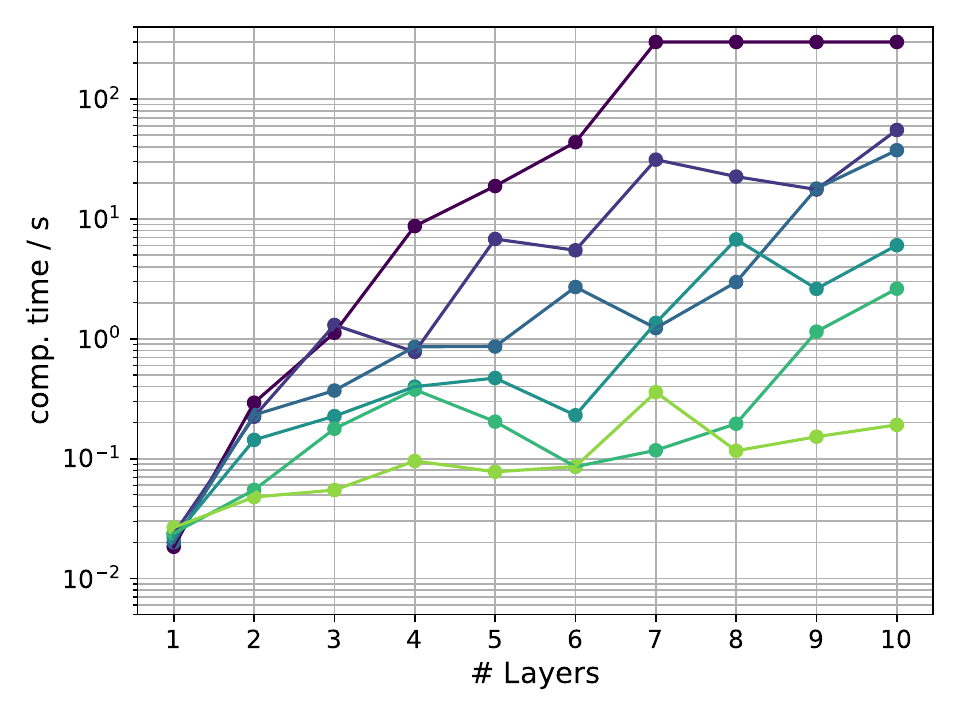}
        \subcaption{Solution time of problem \eqref{prob:minANN} for Himmelblau function.}
    \end{subfigure}

    \begin{subfigure}{\linewidth}
        \centering
        \definecolor{clr1}{HTML}{440154}
\definecolor{clr2}{HTML}{443983}
\definecolor{clr3}{HTML}{31688e}
\definecolor{clr4}{HTML}{21918c}
\definecolor{clr5}{HTML}{35b779}
\definecolor{clr6}{HTML}{90d743}

\begin{tikzpicture} 
    \begin{axis}[%
    hide axis,
    clip bounding box=upper bound,
    xmin=10,
    xmax=50,
    ymin=0,
    ymax=0.4,
    line width= 2,
    scale only axis,
    no markers,
    enlargelimits=false,
    legend columns = 7,
    legend style={draw=none,legend cell align=left}
    ]
    \addlegendimage{empty legend}
    \addlegendentry{$\ell^1$ regularization:}
    \addlegendimage{clr1}
    \addlegendentry{0};
    \addlegendimage{clr2}
    \addlegendentry{1e-7};
    \addlegendimage{clr3}
    \addlegendentry{1e-6};
    \addlegendimage{clr4}
    \addlegendentry{1e-5};
    \addlegendimage{clr5}
    \addlegendentry{1e-4};
    \addlegendimage{clr6}
    \addlegendentry{1e-3};

\end{axis}
\end{tikzpicture}
    \end{subfigure}
    \caption{Mean absolute percentage error on the test set, number of linear regions, percentage of fixed neurons and computation times in problem \eqref{prob:minANN} of trained ANNs with varying number of hidden layers with 25 neurons, trained with different levels of $\ell^1$ regularization.}
    \label{fig:model_statistics_smaller_models}
\end{figure}

Analogous to \Cref{fig:model_statistics}, \Cref{fig:model_statistics_smaller_models} illustrates the results for the ReLU networks with 25 neurons in each hidden layer.

\end{appendices}

\end{document}